\documentclass[11pt,a4paper]{article}
\usepackage{amsbsy,amscd,amsfonts,amssymb,amstext,amsmath,latexsym,theorem}

\usepackage[all]{xy}
\mathchardef\tnode="020E 

\def\arc{
  \hbox{\kern -0.15em
  \vbox{\hrule width 4.5em height 0.6ex depth -0.5 ex}
  \kern -0.33em}}

\def\darc{
  \rlap{\lower0.2ex\arc}{\raise0.2ex\arc}}

\def\stroke#1{
  \kern 0.05em
  \rlap\arc{{\textstyle{#1}}\atop\phantom\arc}
  \kern -0.22em}

\def\dstroke#1{
  \kern 0.05em
  \rlap\darc{{\tiny\textstyle{#1}}\atop\phantom\darc}
  \kern -0.22em}

\def\centerscript#1{
  \setbox0=\hbox{$\tnode$}
  \hbox to \wd0{\hss$\scriptstyle{#1}$\hss}}

\def\node{
  \def\super{}
  \def\sub{}
  \futurelet\next\dolabellednode}

  \let\sp=^
  \let\sb=_

  \def\dolabellednode{%
    \ifx\next\sb\let\next\getsub
    \else
      \ifx\next\sp\let\next\getsuper
      \else\let\next\donode
      \fi
    \fi
    \next}

  \def\getsub_#1{\def\sub{#1}\futurelet\next\dolabellednode}
  \def\getsuper^#1{\def\super{#1}\futurelet\next\dolabellednode}

  \def\donode{%
   \rlap{$\mathop{\phantom\tnode}\limits_{\centerscript{\sub}}^{\centerscript{\super}}$}\tnode}

\def\varcdn{
  \kern -0.03em\vbox{\kern -0.5ex
  \hbox to \wd0{\hss\vrule width 0.04em depth 5.8ex\hss}
  \kern -0.3ex  \hbox{$\tnode$}}}

\newcommand{\Theorem}{Theorem}
\newcommand{\Curtistits}{Curtis-Tits Theorem}
\newcommand{\Titslem}{Tits' Lemma}
\newcommand{\Proposition}{Proposition}
\newcommand{\Lemma}{Lemma}
\newcommand{\Corollary}{Corollary}
\newcommand{\Program}{Program}
\newcommand{\Definition}{Definition}
\newcommand{\Remark}{Remark}
\newcommand{\Example}{Example}
\newcommand{\Examples}{Examples}
\newcommand{\Consequence}{Consequence}
\newcommand{\Fact}{Fact}
\newcommand{\Notation}{Notation}
\newcommand{\Problem}{Problem}
\newcommand{\Conjecture}{Conjecture}
\newcommand{\Technique}{{\sc Technique}}
\theoremstyle{break}
\newtheorem{theorem}{\Theorem}[section]

\newtheorem{proposition}[theorem]{\Proposition}
\newtheorem{lemma}[theorem]{\Lemma}
\newtheorem{corollary}[theorem]{\Corollary}
{\theorembodyfont{\sc}
}

\theoremstyle{plain}
{\theorembodyfont{\rmfamily}

\newtheorem{definition}[theorem]{\Definition}

\newtheorem{example}[theorem]{\Example}

\theoremstyle{break}

}
{\theorembodyfont{\ttfamily}
\theoremstyle{break}

}
\newenvironment{proof}{\noindent {\sl Proof. }}{\hfill \pend \br}

\newcommand{\mc}[1]{\mathcal{#1}}

\newcommand{\pend}{\hfill $\Box$}
\newcommand{\nend}{\hfill $\blacksquare$}

\newcommand{\br}{\vspace{10pt}}

\newcommand{\lb}{\left\{}
\newcommand{\lbr}{\left\lbrack}

\newcommand{\rb}{\right\}}
\newcommand{\rbr}{\right\rbrack}
\renewcommand{\>}{\right>}
\newcommand{\gen}[1]{\left< #1 \right>}

\newcommand{\typ}{{\rm typ}}

\newcommand{\Aut}{{\rm Aut}}

\newcommand{\G}{{\mc{G}}}
\newcommand{\F}{\mathbb{F}}

\newcommand\alp{\alpha}

\newcommand\A{{\cal A}}

\newcommand\B{{\cal B}}

\newcommand\la{\langle}
\newcommand\ra{\rangle}

\newcommand\form{{(\cdot,\cdot)}}

\renewcommand\mod{{\rm mod}\,}
\renewcommand\hat{\widehat}

\newcommand\ie{{i.e.}, }

\newcommand{\SO}{{\mathsf{SO}}}

\newcommand{\sfG}{G}

\begin{document}
\title{{\bf Intransitive geometries} \\ }
\author{Ralf Gramlich \and Hendrik Van Maldeghem}
\date{\today}
\maketitle

\begin{abstract}
A lemma of Tits establishes a connection between the simple connectivity of an incidence geometry and the universal completion
of an amalgam induced by a sufficiently transitive group of automorphisms of that geometry. In the present paper,
we generalize this lemma to intransitive geometries, thus opening the door for numerous applications. We treat
ourselves some amalgams related to intransitive actions of finite orthogonal groups, as a first class of
examples.
\end{abstract}

\section{Introduction}
Amalgams in group theory have proved their importance in the
classification of the finite simple groups (see Sections 28 and 29
of Gorenstein, Lyons, Solomon
\cite{Gorenstein/Lyons/Solomon:1996}). Originally one considers
the amalgam of the maximal parabolic subgroups of a Chevalley
group of rank at least three in its natural action on the associated
building and proves that the universal completion of the amalgam
is (some controlled central extension of) the Chevalley group
itself, see \cite{Curtis:1965}, \cite{Steinberg:1962},
\cite{Tits:1962}, \cite{Tits:1974}. In modern terms, see M\"uhlherr \cite{Muhlherr}, this
essentially is implied by the fact that the building and the opposites geometry of the corresponding twin
building are simply connected.

Since the mid-1970's there has been interest in other types of amalgams as well, see Phan \cite{Phan:1977},
\cite{Phan:1977a}. Somehow miraculously amalgams of (twisted) Chevalley groups over finite fields were studied
that did not come from the action on the building. Aschbacher \cite{Aschbacher} was the first to realize that
Phan's amalgam in \cite{Phan:1977} arises as a version of the amalgam of rank one and rank two parabolics of
the action of $\mathsf{SU}_{n+1}(q^2)$ on the geometry of nondegenerate subspaces of a $(n+1)$-dimensional
unitary vector space over $\F_{q^2}$. In order to prove that the universal completion of the amalgam is the
group under consideration, one complies to a lemma by Tits \cite{Tits:1986}, see also Pasini \cite{Pasini:1985}, saying that this essentially
amounts to checking that the geometry is simply connected and residually connected, under the assumption that
the geometry is flag-transitive.

During the revision of the classification of the finite simple groups there was a demand for a revision of Phan's result \cite{Phan:1977} as well. Das \cite{Das} succeeded
partially and Bennett, Shpectorov \cite{Bennett/Shpectorov}
succeeded completely. After preprints of the latter paper were
circulated around the 2001 conference in honor of Ernie Shult,
things started to develop at a high pace. People finally realized
the connection between M\"uhlherr's \cite{Muhlherr} new proof of
the Curtis-Tits theorem and Aschbacher's \cite{Aschbacher}
geometry for the Phan amalgam. Eventually Hoffman, Shpectorov and
the first author \cite{Gramlich/Hoffman/Shpectorov:2003}
constructed a new geometry resulting in the geometric part of a
completely new Phan-type theorem characterizing central quotients of ${\sf Sp}_{2n}(q)$. Recently the first author
\cite{Gramlich2} provided the group-theoretic part, a
classification of amalgams based on \cite{Bennett/Shpectorov}, thus completing the new Phan-type
theorem. Some remaining open cases over small fields are addressed by Horn in \cite{Horn:2005}, see also Horn, Nickel and the first author \cite{Gramlich/Horn/Nickel}.

Later Bennett joined Hoffman, Shpectorov and the first author
\cite{Bennett/Gramlich/Hoffman/Shpectorov:2003} to develop a
theory for this new sort of geometries, called {\em flipflop
geometries}: Take your favorite spherical building and consider it
as a twin building \`a la Tits \cite{Tits:1992}. The {\em
opposites geometry}, which was used by M\"uhlherr
\cite{Muhlherr} to re-prove the Curtis-Tits theorem, consists
of the pairs of elements of the twin building at codistance one
(the neutral element of the associated Weyl group). A {\em flip}
is an involution of that opposites geometry that interchanges the
positive and the negative part, flips the distances and preserves
the codistance. The flipflop geometry of the opposites geometry
with respect to the flip consists of all those elements of the
opposites geometry that are stabilized (or rather {\em flipped})
by the flip.

In case of Aschbacher's geometry for Phan's theorem the building geometry is the projective space corresponding
to the group $\mathsf{SL}_{n+1}(q^2)$ and the flip is a nondegenerate unitary polarity. The corresponding
flipflop geometry then is the geometry on the nondegenerate subspaces of the projective space with respect to
the polarity. Indeed, being opposite means that a subspace and its polar have empty intersection which in turn
means that the subspace in question is nondegenerate.

The rank of this geometry is always higher than the one of the
associated building, and hence this approach covers more groups.
This idea works fine for the unitary groups (see Aschbacher
\cite{Aschbacher}, Das \cite{Das}, Bennett, Shpectorov
\cite{Bennett/Shpectorov}) and for the symplectic groups (see Das
\cite{Das:1998} (finite fields, odd characteristic), Das
\cite{Das:2000} (finite fields, even characteristic), Hoffman,
Shpectorov and the first author
\cite{Gramlich/Hoffman/Shpectorov:2003} (finite fields of size at
least $8$; a by-product of the new geometry), and the first author
\cite{Gramlich3} (all fields)) although, strictly speaking, the
symplectic forms do not yield a flipflop geometry. However, for
the orthogonal ones over finite fields, we run into problems since
the geometry of nondegenerate spaces is, in general, not
flag-transitive. The flag-transitive case for forms of Witt index
at least one, i.e., over quadratically closed fields has been
settled by Altmann \cite{Altmann:2003}. See also Altmann and the
first author \cite{Altmann/Gramlich} for the same results and some
extensions to real closed fields.

As said before, in order to prove that the universal completion of the amalgam is the group under
consideration, one complies to a lemma by Tits \cite{Tits:1986} saying that this essentially amounts to
checking that the geometry is simply connected and residually connected, under the assumption of
flag-transitivity. For intransitive geometries one can try to find a flag-transitive subgeometry and to prove
that this subgeometry is simply connected and residually connected. However, flag-transitive subgeometries of
the geometry of nondegenerate subspaces of a finite orthogonal classical group are not known to be simply
connected in low rank. Adam Roberts has produced some results in this direction for rank (of the geometry) at least four, see \cite{Rob}.

Hence, to overcome the difficulties that occur in rank three, one should generalize the theory of amalgams either to non-flag-transitive geometries, or to non-simply connected ones. Since the former is more direct (the latter
requires suitable flag-transitive subgeometries and involves nontrivial covers of these subgeometries; this
seems to be less direct), we have chosen to try that. The key idea is to use a theorem by Stroppel
\cite{Stroppel:1993}, which seems not to be so well known, but is very useful in this context. We also discuss
the more difficult and more general problem of the amalgam of rank $k$ parabolics in non-flag-transitive
geometries. It actually turns out that the most natural results occur if one abandons thinking in amalgams of
rank $k$ parabolics, but adopts thinking in amalgams of certain {\em shapes} instead. We then apply our theory
to the orthogonal classical groups and give many examples.

In an appendix we apply our theory reporting on recent research by Hoffman and Shpectorov
\cite{Hoffman/Shpectorov} on an interesting amalgam for $\mathsf{G_2}(3)$ coming from an intransitive geometry related
to the sporadic simple Thompson group. Our approach shows how powerful the established theory of the present
paper is.

We conclude this introduction by the remark that in the
mid-1980's, using functional analysis and Lie theory, Borovoi
\cite{Borovoi:1984} and Satarov \cite{Satarov:1985} have obtained
related universal completion results for amalgams in compact Lie
groups. In this case, however, the geometry acted on is the
building, so their results on compact Lie groups follow
immediately from the simple connectivity of the building. The
classification strategy for amalgams from
\cite{Bennett/Shpectorov} and \cite{Gramlich2} was used by the
first author in \cite{Gramlich4} when providing a classification
of the amalgams from \cite{Borovoi:1984} and \cite{Satarov:1985},
yielding a Phan-type theorem for compact Lie groups. Gl\"ockner and the first author \cite{Gloeckner/Gramlich} proved recently that these amalgams in fact suffice to reconstruct the Lie group topology as well.

\medskip \noindent
{\bf Acknowledgement:} The authors would like to express their gratitude to Antonio Pasini for his detailed
remarks and suggestions, requiring corrections and improvements at various places. The authors are also grateful to Max Horn for additional remarks and comments.

\section{Preliminaries}

In this section, we define the notions and review the results that we will need to develop our theory. This
section has been inspired by \cite{Buekenhout/Cohen}, \cite{Seifert/Threlfall:1934}, \cite{Serre:2003}.

\subsection{Coset pregeometries}

\begin{definition}[Pregeometry, geometry]
 A \textbf{pregeometry} $\mc{G}$ \textbf{over the set $I$} is a
triple $(X,*,\typ)$ consisting of a set $X$, a symmetric and
reflexive \textbf{incidence relation} $*$, and a surjective
\textbf{type function} $\typ:X\rightarrow I$, subject to the
following condition:
\begin{description}\item[(Pre)] If $x*y$ with $\typ(x)=\typ(y)$, then
$x=y$.
\end{description} The set $I$ is usually called the \textbf{type set}.
A \textbf{flag} in $X$ is a set of pairwise incident elements. The
\textbf{type} of a flag $F$ is the set $\typ(F):=\{\typ(x):x\in
F\}$. A \textbf{chamber} is a flag of type $I$, a {\bf pennant} is a flag of cardinality three. The \textbf{rank}
of a flag $F$ is $|\typ(F)|$ and the \textbf{corank} is equal to
$|I\setminus\typ(F)|$.

A \textbf{geometry} is a pregeometry with
the additional property that
\begin{description}
\item[(Geo)] every flag is contained in a chamber.
\end{description}

The pregeometry $\mc{G}$ is \textbf{connected} if the graph
$(X,*)$ is connected.
\end{definition}

\begin{definition}[Lounge, hall]
Let $\mathcal{G} = (X,*,\typ)$ be a pregeometry over $I$. A subset $W$ of $X$ is called a {\bf lounge} if each
subset $V$ of $W$ for which $\typ : V \rightarrow I$ is a injection, is a flag. A lounge $W$ with $\typ(W) = I$
is called a {\bf hall}.
\end{definition}

\begin{definition}[Residue]
Let $F$ be a flag of $\mc{G}$, let us say of type $J\subseteq I$. Let $J'$ be the set of all $j'\in I\setminus
J$ such that there exists $x\in X$ with $\typ(x)=j'$ and with $F\cup\{x\}$ a flag. Also, let $X'$ be the set of
all $x\in X$ with $\typ(x)\in J'$ such that $F\cup\{x\}$ is a flag. Then the \textbf{residue $\mc{G}_F$ of $F$}
is the pregeometry
$$(X',*_{|X'\times X'},\typ_{|X'})$$ over $J'$.
\end{definition}

\begin{definition}[Automorphism]
Let $\mc{G}=(X,*,\typ)$ be a pregeometry over $I$. An
\textbf{automorphism} of $\mc{G}$ is a permutation $\sigma$ of $X$
with $\typ(\sigma(x))=\typ(x)$, for all $x\in X$, and with
$\sigma(x)*\sigma(y)$ if and only if $x*y$, for all $x,y\in X$.
\begin{tabbing}
A
group $G$ of automorphisms of $\mc{G}$ is called \quad \= if for each pair of
flags $c$, $d$ with \kill
Moreover, a
group $G$ of automorphisms \> \\
of $\mc{G}$ is called \> \\ \> if for each pair of flags $c$, $d$ with \\
 {\bf flag-transitive}, \> $\typ(c) = \typ(d)$, \\
 {\bf chamber-transitive}, \> $\typ(c) = I = \typ(d)$, \\
{\bf pennant-transitive}, \> $|\typ(c)| = 3 = |\typ(d)|$ and $\typ(c) = \typ(d)$, \\
{\bf incidence-transitive}, \> $|\typ(c)| = 2 = |\typ(d)|$ and $\typ(c) = \typ(d)$, \\
{\bf vertex-transitive} \> $|\typ(c)| = 1 = |\typ(d)|$ and $\typ(c) = \typ(d)$ \\ \\
there exists a $\sigma \in G$ with $\sigma(c) = d$.\>
\end{tabbing}

If the group of {all} automorphisms of $\mc{G}$ is flag-transitive, chamber-transitive, incidence-transitive or
vertex-transitive, then we say that $\mc{G}$ is {\bf flag-transitive}, {\bf chamber-transitive}, {\bf
incidence-transitive} or {\bf vertex-transitive}, respectively.
\end{definition}

The emphasis of the present paper is on geometries that are not
vertex-transitive, and which we will call \textbf{intransitive}.
Therefore, we first have a look how one can describe such a
geometry group-theoretically.

\begin{definition}[Coset Pregeometry] \label{cospre}
Let $I$ be a set and let $(T_i)_{i \in I}$ be a family of mutually disjoint sets. Also, let $G$ be a group and
let $(G^{t,i})_{t \in T_i, i \in I}$ be a family of subgroups of $G$. Then
$$\left( \left\{(C,t):t\in T_i\mbox{ for some }i\in I, C\in G/G^{t,i} \right\},*,\typ \right)$$ with $\typ (C,t) = i$ if $t\in T_i$, and
\begin{description}
\item[(Cos)] $gG^{t,i} *
hG^{s,j}$ if and only if $gG^{t,i} \cap hG^{s,j} \neq \emptyset$
and either $i \neq j$ or $(t,i) = (s,j)$
\end{description}
is a pregeometry over $I$,
the {\bf coset pregeometry of $G$} with respect to $(G^{t,i})_{t
\in T_i, i \in I}$. Since the type function is completely
determined by the indices, we also denote the coset pregeometry of
$G$ with respect to $(G^{t,i})_{t \in T_i, i \in I}$ by
$$((G/G^{t,i}\times\{t\})_{t \in T_i, i \in I}, *).$$
The family $(G^{t,i},t)_{t \in T_i, i \in I}$ forms a lounge, even a hall. If $|T_i| = 1$ for all $i \in I$,
then we write $G_i$ instead of $G^{t,i}$, and we put canonically $T_i=\{i\}$. The family $(G_i,i)_{i \in I}$
forms a chamber of the coset geometry, called the {\bf base chamber}. For $J\subseteq I$, we also write
$G_J=\bigcap_{j\in J}G_j$.
\end{definition}

Certainly, any coset pregeometry with $|T_i| = 1$ for all $i \in I$, which means nothing else than being
vertex-transitive, is incidence-transitive. Indeed, if $gG_{i} \cap hG_{j} \neq \emptyset$, then choose $a \in
gG_i \cap hG_j$. It follows $aG_i = gG_i$ and $aG_j = hG_j$ and therefore the automorphism $a^{-1}$ maps the
incident pair $gG_{i}$, $hG_{j}$ onto the incident pair $G_{i}$, $G_{j}$.

\medskip
Note that the residue of a coset pregeometry in general is not a coset pregeometry. See \cite{Buekenhout/Cohen} for a number of conditions under which it in fact \emph{is} a coset pregeometry.

\medskip
Similar to the characterizations of vertex-transitivity there
exist a large number of group-theoretic characterizations of
various geometric properties of coset geometries, see e.g.\
\cite{Buekenhout/Cohen}. The following one, the characterization
of connectivity, is an easy but crucial observation for studying
amalgams.

\begin{theorem}[inspired by Buekenhout/Cohen \cite{Buekenhout/Cohen}] \label{char conn}
Let $I \neq \emptyset$. The coset pregeometry $((G/G^{t,i}\times \{t\})_{t \in T_i, i \in I}, *)$ is connected if
and only if
$$G= \langle G^{t,i} \mid  i \in I, t \in T_i \rangle.$$
\end{theorem}

\begin{proof}
Suppose that $\G$ is connected. Take $i\in I$ and $t \in T_i$.
If $a\in G$, then there is a path
$$1G^{t,i}, a_0 G^{t_0, i_{0}},
a_1 G^{t_1,i_{1}},
a_{2} G^{t_2,i_{2}}, \ldots,
a_{m} G^{t_m,i_{m} }, a G^{t,i}$$
connecting the elements $1G^{t,i}$ and
$a G^{t,i}$ of $\G$. Now
$$a_k G^{t_k,i_k} \cap a_{k+1} G^{t_{k+1},i_{k+1}} \neq \emptyset,$$
so $$a_k^{-1} a_{k+1}\in G^{t_k,i_k} G^{t_{k+1},i_{k+1}}$$
for $k=0 ,\ldots, m-1$.
Hence
$$a = (1^{-1}a_0)( a_0^{-1} a_1 ) \cdots
(a_{m-1}^{-1} a_{m})(a_m^{-1}a) \in  G^{t,i}G^{t_0,i_0} \cdots G^{t_{m-1},i_{m-1}} G^{t_m, i_{m}}G^{t,i},$$ and
so $a\in \langle G^{t,i} \mid  i \in I, j \in T_i \rangle$. The converse is obtained by reversing the above
argument. The only difficulty that can occur is that $g_1G^{t_1,i}$ and $g_2G^{t_2,i}$ are not incident, even
if $g_1G^{t_1,i} \cap g_2G^{t_2,i} \neq \emptyset$. This can be remedied by including some coset $gG^{t,j}$, $j
\neq i$, between $g_1G^{t_1,i}$ and $g_2G^{t_2,i}$ into the chain of incidences, where $g \in g_1G^{t_1,i} \cap
g_2G^{t_2,i}$.
\end{proof}

Now we turn to the question which pregeometries actually are coset pregeometries. Stroppel gave the
answer in \cite{Stroppel:1993}. To this end let us introduce
the notion of the sketch of a pregeometry.

\begin{definition}[Sketch] \label{sketch}
Let $\mc{G} = (X,*,\typ)$ be a pregeometry over $I$, let $G$ be a
group of automorphisms of $\mc{G}$, and let $W \subset X$ be a set
of $G$-orbit representatives of $X$. We write $$W=\bigcup_{i\in
I}W_i$$ with $W_i\subseteq\typ^{-1}(i)$. The {\bf sketch of $\G$
with respect to $G$ and $W$} is the coset geometry
$$((G/G_w\times\{w\})_{w\in W_i, i\in I},*).$$
\end{definition}

Recall that two actions $$\phi
: G \to \Aut\ M \quad \mbox{ and } \quad \phi' : G \to \Aut\ M'$$ are said to be
{\bf equivalent}
if there is an isomorphism
$\psi : M \to M'$ such that
$\psi \circ \phi(g) \circ \psi^{-1} = \phi'(g)$ for each $g\in G$ or,  equivalently,
$\psi \circ \phi(g) = \phi'(g) \circ \psi$ for all $g\in G$. In this case,
we shall also say that $M$ and $M'$ are {\bf isomorphic $G$-sets}.

\begin{theorem}[Stroppel's reconstruction theorem \cite{Stroppel:1993}] \label{isomorphism theorem2}
Let $\mc{G} = (X,*,\typ)$ be a pregeometry over $I$ and let $G$ be a group of automorphisms of $\G$. For each
$i \in I$ let $$w^i_1, \ldots, w^i_{t_i}$$ be $G$-orbit representatives of the elements of type $i$ of $\G$
such that
\begin{enumerate}
\item $W := \bigcup_{i \in I} \lb w_1^i, \ldots, w_{t_i}^i \rb$ is a hall and, 
\item if $V \subseteq W$ is a
chamber, the action of $G$ on the pregeometry over $\typ(V)$ consisting of all elements of the $G$-orbits $G.x$,
$x \in V$, is incidence-transitive.
\end{enumerate}
Then the bijection $\Phi$ between the sketch of $\G$ with respect to $G$ and $W$ and the pregeometry $\G$ given
by $$gG_{w_j^i} \mapsto gw_j^i$$ is an isomorphism between pregeometries and an isomorphism between $G$-sets.
\pend
\end{theorem}

For a vertex-transitive group $G$, the previous theorem is just
the isomorphism theorem of incidence-transitive pregeometries, see
\cite{Buekenhout/Cohen}.

\medskip
The geometry consisting of the $G$-orbits $G.x$ of elements of some fixed chamber $V \subseteq W$ as in
(ii) of the theorem is called the {\bf orbit geometry for $(\G,G,V)$}.

\subsection{Fundamental group and simple connectivity}

\begin{definition}[Fundamental group] Let $\mc{G}$ be a connected pregeometry.
A path of length $k$ in the geometry is a sequence of elements
$(x_0, \ldots, x_k)$ such that $x_i$ and $x_{i+1}$ are incident, $0 \leq i \leq k-1$. A {\bf cycle} based at an
element $x$ is a path in which $x_0 = x_k = x$. Two paths based at the same vertex are {\bf
homotopically equivalent} if one can be obtained from the other
via the following operations (called {\bf elementary homotopies}):
\begin{enumerate}
\item inserting or deleting a repetition (\ie a cycle of length 1),
\item inserting or deleting a return (\ie a cycle of length 2), or
\item inserting or deleting  a
triangle (\ie a cycle of length 3).
\end{enumerate}
  The equivalence classes of
cycles based at an element $x$ form a group under the operation
induced by concatenation of cycles. This group is called the {\bf
fundamental group} of $\G$ and denoted by $\pi_1(\G, x)$.

A cycle based at $x$ that is homotopically equivalent to the
trivial cycle $(x)$ is called \textbf{null-homotopic}. Every cycle
of length $1$, $2$, or $3$ is null-homotopic.
\end{definition}

\begin{definition}[Covering]
Suppose $\G$ and $\hat\G$ are two connected geometries over the
same type set and suppose $\phi:\hat\G\rightarrow\G$ is a {\bf
homomorphism} of geometries, \ie $\phi$ preserves the types and
sends incident elements to incident elements.  A surjective
homomorphism $\phi$ between connected geometries $\hat \G$ and
$\G$ is called a {\bf covering} if and only if for every nonempty
flag $\hat F$ in $\hat\G$ the mapping $\phi$ induces an
isomorphism between the residue of $\hat F$ in $\hat\G$ and the
residue of $F=\phi(\hat F)$ in $\G$.  Coverings of a geometry
correspond to the usual topological coverings of the flag complex.
It is well-known and easy to see that a surjective homomorphism
$\phi$ between connected geometries $\hat\G$ and $\G$ is a
covering if and only if for every element $\hat x$ in $\hat \G$
the map $\phi$ induces an isomorphism between the residue of $\hat
x$ in $\hat\G$ and the residue of $x = \phi(\hat x)$ in $\G$. If
$\phi$ is an isomorphism, then the covering is said to be
\textbf{trivial}.
\end{definition}

Consider the geometry via its colored incidence graph and recall the following results from the theory of
simplicial complexes.

\begin{theorem}[Chapter 8 of Seifert/Threlfall \cite{Seifert/Threlfall:1934}]
Let $\G$ be a connected geometry and let $x$ be an element of
$\G$. The group $\pi_1(\G, x)$ is trivial if and only if all coverings of $\G$ are trivial. \pend
\end{theorem}

A geometry satisfying the equivalent conditions in the previous
theorem is called \textbf{simply connected}.

\medskip The following construction can also be found in Chapter 8 of \cite{Seifert/Threlfall:1934}.

\begin{definition}[Fundamental cover]
Let $\Gamma$ be a connected graph and let $x$ be some vertex of $\Gamma$. The {\bf fundamental cover
$\widehat{\Gamma}$ of $\Gamma$ based at $x$} is defined as follows: The vertices of $\widehat{\Gamma}$ are the
homotopy classes of paths of $\Gamma$ based at $x$ where two vertices $\lbr \gamma_1 \rbr$ and $\lbr \gamma_2
\rbr$ of $\widehat{\Gamma}$ are adjacent if and only if $\lbr \gamma_1^{-1} \gamma_2 \rbr = \lbr t_1 t_2 \rbr$
where $t_1$ is the terminal vertex of $\gamma_1$, $t_2$ is the terminal vertex of $\gamma_2$, and $t_1$ and
$t_2$ are adjacent in $\Gamma$.
\end{definition}

\begin{definition}[Universal covering] \label{univcov}
Let $\Gamma$ and $\widehat\Gamma$ be connected graphs and let $x \in \Gamma$, $\widehat{x} \in \widehat\Gamma$
be vertices. A covering $$\pi : \hat\Gamma \rightarrow \Gamma$$ mapping $\widehat{x}$ onto $x$ is called {\bf
universal} if, for any covering $$\alpha : \Gamma_1 \rightarrow \Gamma \quad \mbox{ and any } \quad x_1 \in
\alpha^{-1}(x),$$ there exists a unique covering map $$\beta : \hat\Gamma \rightarrow \Gamma_1$$ with $\pi =
\alpha \circ \beta$ and $\beta(\widehat{x}) = x_1$.
\end{definition}

$$\xymatrix{
(\hat{\Gamma},\widehat{x}) \ar[r]^\beta \ar[dr]_\pi & (\Gamma_1,x_1) \ar[d]^\alpha \\
& (\Gamma,x)
}$$

\begin{theorem}[Chapter 8 of Seifert/Threlfall \cite{Seifert/Threlfall:1934}] \label{1universal property}
Let $\Gamma$ be a connected graph, let $x$ be a vertex of $\Gamma$, and let $\widehat\Gamma$ be the fundamental
cover of $\Gamma$ based at $x$. Then the fundamental covering $\pi : \hat\Gamma \rightarrow \Gamma$ is
universal. \pend
\end{theorem}

\subsection{Amalgams}

\begin{definition}[Amalgam]
An {\bf amalgam of groups} $\A$ over a finite set $I = \lb 0, 1,
\ldots, n \rb$ and associated nonempty sets $J_i$, $i \in I$, is a
family of groups $(G_{j,i})_{j \in J_i, i \in I}$ with
monomorphisms, called {\bf identifications},
$$\phi_{j_i,i}^{j_{i+1},i+1} : G_{j_i,i} \rightarrow
G_{j_{i+1},i+1}$$ for some $(j_i,i)$ and $(j_{i+1},i+1)$ such that
for each $G_{j_i,i}$ there exist identifications whose composition
embeds $G_{j_i,i}$ into some $G_{j_n,n}$.
\end{definition}

\begin{example} \label{2.18}
An amalgam with $I = \lb 0, 1, 2 \rb$, $J_{0} = \lb 1, 2 \rb$, $J_{1} = \lb 1, 2, 3, 4 \rb$, $J_{2} = \lb 1, 2,
3, 4 \rb$ can be depicted in the following diagram. The identification maps are given by arrows.
$$
\xymatrix{
& G_{1,1} \ar[r]^{\phi_{1,1}^{1,2}} \ar[ddr] & G_{1,2} \\
& G_{2,1} \ar[r] \ar[dr]  & G_{2,2} \\
G_{1,0} \ar[uur]^{\phi_{1,0}^{1,1}} \ar[ddr] & & G_{3,2} \\
G_{2,0} \ar[uur]_>>>>>>>>>>{\phi_{2,0}^{2,1}} \ar[r] \ar[dr] & G_{3,1} \ar[uur] \ar[r] & G_{4,2} \\
& G_{4,1} \ar[ur]_{\phi_{4,1}^{4,2}} \ar[uur]_>>>>>{\phi_{4,1}^{3,2}} &
}
$$
Note that the definition of an amalgam does \emph{not} imply $$\phi_{2,1}^{3,2} \circ \phi_{2,0}^{2,1} =
\phi_{4,1}^{3,2} \circ \phi_{2,0}^{4,1}$$ in the above example.

\medskip
Two amalgams $\A$ and $\B$ are {\bf similar} if they share the same set $I$, the same sets $J_i$ and if for all
$(j_i,i)$ and $(j_{i+1},i+1)$ the identification $_\A\phi_{j_i,i}^{j_{i+1},i+1}$ exists if and only if the
identification $_\B\phi_{j_i,i}^{j_{i+1},i+1}$ exists, i.e., if they can be depicted by the same diagram.
\end{example}

\begin{definition}[Homomorphism]
Let $\A=(G_{j,i})_{j,i}$ and $\B=(H_{j,i})_{j,i}$ be similar amalgams.  A map $\psi : \sqcup \A \rightarrow
\sqcup \B$ will be called an {\bf amalgam homomorphism from $\A$ to $\B$} if
\begin{enumerate}
\item for every $i\in I$ and $j \in J_i$ the restriction of $\psi$ to $G_{j,i}$ is a homomorphism from
$G_{j,i}$ to $H_{j,i}$ and \item $\psi \circ {}_\A\phi_{j_i,i}^{j_{i+1},i+1} = {}_\B\phi_{j_i,i}^{j_{i+1},i+1}
\circ \psi_{|G_{j_i,i}}$ in case the respective identifications exist.
\end{enumerate}
If $\psi$ is bijective and its inverse map $\psi^{-1}$ is also an amalgam homomorphism, then $\psi$ is called an {\bf
amalgam isomorphism}.\index{isomorphism}  An {\bf automorphism}\index{automorphism} of $\A$ is an isomorphism
of $\A$ onto itself.  As usual, the automorphisms of $\A$ form the
{\bf automorphism group}, $\Aut(\A)$.
\end{definition}

\begin{definition}[Quotient, cover]
An amalgam $\B=(H_{j,i})_{j,i}$ is a {\bf quotient}\index{quotient} of the amalgam $\A=(G_{j,i})_{j,i}$ if
there is an amalgam homomorphism $\pi$ from $\A$ to $\B$ such that the restriction of $\pi$ to any $G_{j,n}$
maps $G_{j,n}$ onto $H_{j,n}$.  The map $\pi : \sqcup \mc{A} \rightarrow \sqcup \mc{B}$ is called a {\bf
covering}, $\mc{A}$ is called a {\bf cover}\index{cover}\index{amalgam!cover} of $\mc{B}$. Two coverings
$(\A_1,\pi_1)$ and $(\A_2,\pi_2)$ of $\A$ are called {\bf
equivalent}\index{amalgam!cover!equivalent}\index{cover!equivalent} if there is an isomorphism $\psi$ of $\A_1$
onto $\A_2$, such that $\pi_1=\pi_2 \circ \psi$.
\end{definition}

Notice that a covering $\pi : \sqcup \mc{A} \rightarrow \sqcup \mc{B}$ between amalgams need not map $G_{j,i}$
surjectively onto $H_{j,i}$ for $i \neq n$.

\begin{definition}[Completion]\label{completiondef}
Let $\A$ be an amalgam. A pair $(G,\pi)$ consisting of a group $G$ and a {map} $\pi: \sqcup\A\rightarrow G$ is
called a {\bf completion} of $\A$, and $\pi$ is called a {\bf completion map}, if
\begin{enumerate}
\item for all $i\in I$ and $j \in J_i$  the restriction of $\pi$ to $G_{j,i}$ is a homomorphism of $G_{j,i}$ to
$G$; \item $\pi_{|G_{j_{i+1},i+1}}\circ \phi_{j_i,i}^{j_{i+1},i+1} = \pi_{|G_{j_{i},i}}$ if the corresponding
identification exist; and \item $\pi(\sqcup \A)$ generates $G$.
\end{enumerate}
A completion is called {\bf faithful} if for each $i \in I$ and $j \in J_i$ the restriction of $\pi$ to
$G_{j,i}$ is injective.
\end{definition}

Coming back to Example~\ref{2.18}, the definition of a completion {does} require that $$\pi_{|G_{3,2}} \circ
\phi_{2,1}^{3,2} \circ \phi_{2,0}^{2,1} = \pi_{|G_{3,2}} \circ \phi_{4,1}^{3,2} \circ \phi_{2,0}^{4,1},$$
although by definition of an amalgam we do {not necessarily} have $$\phi_{2,1}^{3,2} \circ \phi_{2,0}^{2,1} =
\phi_{4,1}^{3,2} \circ \phi_{2,0}^{4,1} .$$

\begin{proposition} \label{universal}
Let $\A=(G_{j,i})_{j,i}$ be an amalgam of groups, let $F(\A) = \gen{(u_g)_{g\in\A}}$ be the free group on the
elements of $\A$ and let $$S_1 = \lb u_xu_y=u_z,\mbox{ whenever $xy=z$ in some $G_{j,i}$} \rb$$ and $$S_2 = \lb
u_x = u_y,\mbox{ whenever $\phi(x)=y$ for some identification $\phi$} \rb $$ be relations for $F$. Then for
each completion $(G,\pi)$ of $\A$ there exists a unique group epimorphism $$\hat \pi : \mc{U}(\mc{A})
\rightarrow G$$ with  $\pi = \hat \pi \circ \psi$ where $$\mc{U}(\A) = \gen{(u_g)_{g\in\A} \mid S_1,S_2} \mbox{
and } \psi : \sqcup\mc{A} \rightarrow \mc{U}(\A) : g \mapsto u_g.$$
$$\xymatrix{
\sqcup\A \ar[r]^\psi \ar[dr]_\pi & \mc{U}(\mc{A}) \ar[d]^{\hat\pi} \\
&  G
}$$
\end{proposition}

\begin{proof}
The map $\A$ to $\mc{U}(\A)$ given by $\psi: g\mapsto u_g$ turns the group $\mc{U}(\A)$ into a completion of
$\A$.  If $(G,\pi)$ is an arbitrary completion of $\A$ then the map $$\hat \pi : u_g\mapsto\pi(g)$$ leads to a
group epimorphism $\hat\pi$ from $\mc{U}(\A)$ to $G$ because $$\hat \pi(u_g u_h) = \hat \pi(u_{gh}) = \pi(gh) =
\pi(g) \pi(h) = \hat \pi(u_g) \hat \pi(u_h)$$ if $u_{gh}$ exists; otherwise define $$\hat \pi(u_g u_h):=\pi(g)
\pi(h) = \hat \pi(u_g) \hat \pi(u_h).$$ Clearly, $\hat \pi$ is uniquely determined by the requirement $$\pi(g)
= (\hat \pi \circ \psi)(g) = \hat \pi(u_g).$$
\end{proof}

\begin{definition}[Universal Completion]
Let $\A=(G_{j,i})_{j,i}$ be an amalgam of groups. Then $$\psi : \sqcup\mc{A} \rightarrow \mc{U}(\A) : g \mapsto
u_g$$ for $\mc{U}(\A)$ as in Proposition \ref{universal} is called the {\bf universal completion of $\A$}. The
amalgam $\A$ {\bf collapses} if $\mc{U}(\A)=1$
\end{definition}

Notice that if $\B$ is a quotient of $\A$ then $\mc{U}(\B)$ is isomorphic to a factor group of $\mc{U}(\A)$.
In particular, if $\B$ does not collapse then neither does $\A$. Also, an amalgam $\mc{A}$ admits a faithful
completion if and only if its universal completion is faithful.

\begin{definition}[Amalgams for transitive geometries] \label{defrank}
Suppose $\G$ is a geometry and $G\le\Aut\ \G$ is an incidence-transitive group.  Corresponding to $\G$ and $G$
and some chamber $F$, there is an amalgam $\A=\A(\G,G,F)$, the {\bf amalgam of parabolics with respect to $\G$,
$G$, $F$}, defined as the family $(G_E)_{\emptyset\neq E \subseteq F}$, where $G_E$ denotes the stabilizer of
$E \subseteq F$ in $G$, together with the natural inclusions as identification maps.  In case $G$ is
flag-transitive, the amalgam $\A$ is independent (up to conjugation) of the choice of $F$.

The group $G_E$ is called a \textbf{parabolic subgroup of rank $|F\setminus E|$}.

If $|I| = n$ is finite and $k < n$ the amalgam $\A_{(k)} = \A_{(k)}(\G,G,F)$ is the {subamalgam} of $\A$
consisting of all parabolics of rank less or equal $k$. It is called the {\bf amalgam of rank $k$ parabolics}.
Of course, $\A_{(n-1)} = \A$.

More generally, for $F$ as above suppose $\mc{W} \subseteq 2^F$ such that $2^F \ni U' \supseteq U \in \mc{W}$
implies $U' \in \mc{W}$, i.e., $\mc{W}$ is a subset of the power set of $F$ that is closed under passing to
supersets. A set $\mc{W} \subseteq 2^F$ with those properties is called a {\bf shape}. The {\bf amalgam of shape
$\mc{W}$ with respect to $\mc{G}$, $G$, $F$} is the family $(G_U)_{\emptyset\neq U \in \mc{W}}$, where $G_U$ is
the stabilizer of $U \in \mathcal{W}$ in $G$, with the natural inclusion maps as identification maps. It is
denoted by $\mc{A}_\mc{W}(\mc{G},G,F)$.
\end{definition}

For example, note that the amalgam $\mc{A} = \mc{A}(\mc{G},G,F)$ is an amalgam of shape
$2^F\setminus\{\emptyset\}$.

\begin{definition}[Amalgams for intransitive geometries] \label{defrank2}
Suppose $\mc{G} = (X,*,\typ)$ is a geometry over $I$, the group $G$ is a group of automorphisms of $\G$, and
for each $i \in I$ let $w^i_1$, $\ldots$, $w^i_{t_i}$ be $G$-orbit representatives of the elements of type $i$
of $\G$ such that
\begin{enumerate}
\item $W := \bigcup_{i \in I} \lb w_1^i, \ldots, w_{t_i}^i \rb$ is a hall and, 
\item if $V \subseteq W$ is a
chamber, the action of $G$ on the pregeometry over $I$ consisting of all elements of the $G$-orbits $G.x$, $x \in
V$, is incidence-transitive.
\end{enumerate}
Then the amalgam $\mc{A} = \mc{A}(\mc{G},G,W)$ is defined as the family $(G_U)_{\emptyset\neq U \subseteq W
\mbox{ a flag}}$, where $G_U$ denotes the stabilizer of $U \subseteq W$ in $G$ with the natural inclusion maps
as identification maps.

For example, let $\G$ be a rank three geometry with $W$ equal to $p$, $q$, $l$, $\pi$, where $p$ and $q$ are elements of the same type and $p$, $l$, $\pi$ each have different types. Then the amalgam of
parabolics looks as follows:
$$
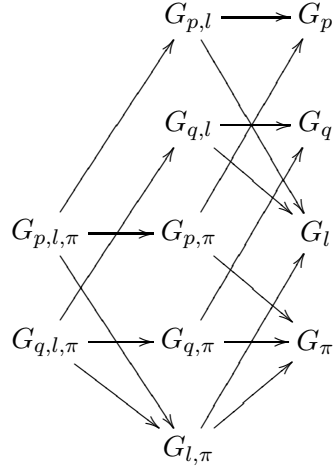
\xymatrix{
& G_{p,l} \ar[r] \ar[ddr] & G_p \\
& G_{q,l} \ar[r] \ar[dr]  & G_q \\
G_{p,l, \pi} \ar[uur] \ar[r] \ar[ddr] & G_{p,\pi} \ar[uur] \ar[dr] & G_l \\
G_{q,l,\pi} \ar[uur] \ar[r] \ar[dr] & G_{q,\pi} \ar[uur] \ar[r] & G_\pi \\
& G_{l,\pi} \ar[ur] \ar[uur] &
}
$$

If $|I| = n$ is finite and $k < n$ the amalgam $\A_{(k)} = \A_{(k)}(\G,G,W)$ is the subamalgam of $\A$
consisting of all parabolics of rank less or equal $k$. It is called the {\bf amalgam of rank $k$ parabolics}.
Of course, $\A_{(n-1)} = \A$.

More generally, for $W$ as above suppose $\mc{W} \subseteq 2^W$ with the properties that each $U \in \mc{W}$ is
a flag and if $U' \subseteq W$ is a flag with $U' \supseteq U \in \mc{W}$, then also $U' \in \mc{W}$, i.e.,
$\mc{W}$ is a subset of the power set of $W$ consisting of flags that is closed under passing to
``superflags''. A set $\mc{W} \subseteq 2^W$ with those properties is called a {\bf shape}. The {\bf amalgam of
shape $\mc{W}$ for $(\mc{G},G,W)$} is defined on the family $(G_U)_{U \in \mc{W}}$ with the natural inclusion
maps as identification maps.  It is denoted by $\mc{A}_\mc{W}(\mc{G},G,W)$.
\end{definition}

For example, note that the amalgam $\mc{A} = \mc{A}(\mc{G},G,W)$ is an amalgam of shape
$$\{ F \in 2^W \mid F \mbox{ nonempty flag} \}.$$

\section{Theory of intransitive geometries}

We now use the foregoing notions, definitions and basic results to develop some theory of intransitive
 geometries, that results in criteria to conclude that certain completions of certain amalgams are
universal.

\begin{theorem}[Fundamental theorem of geometric covering theory] \label{mainthmcover}
Let $\mc{G} = (X,*,\typ)$ be a connected geometry over $I$ of rank at least three and let $G$ be a group of
automorphisms of $\G$. For each $i \in I$ let $$w^i_1, \ldots, w^i_{t_i}$$ be $G$-orbit representatives of the
elements of type $i$ of $\G$ such that
\begin{enumerate}
\item $W := \bigcup_{i \in I} \lb w_1^i, \ldots, w_{t_i}^i \rb$ is a hall and, 
\item if $V \subseteq W$ is a
chamber, the action of $G$ on the pregeometry over $\typ(V)$ consisting of all elements of the $G$-orbits $G.x$,
$x \in V$, is incidence-transitive and pennant-transitive.
\end{enumerate}
Let $\mc{A} = \mc{A}(\G,G,W)$ be the amalgam of parabolics. Then the coset pregeometry $$\hat\G =
((\mc{U}(\mc{A})/G_{w^i_j}\times\{w^i_j\})_{1 \leq j \leq t, i \in I},*)$$ is a simply connected geometry that
admits a universal covering $\pi : \hat\G \rightarrow \G$ induced by the natural epimorphism $\mc{U}(\mc{A})
\rightarrow G$. Moreover, $\mc{U}(\mc{A})$ is of the form $\pi_1(\G).G$.
\end{theorem}

\begin{proof}
First notice that, since $\mc{G}$ is connected, $G$ is generated by all its parabolics (different from $G$) by
Theorem~\ref{char conn}. Hence by Definition~\ref{completiondef}, $G$ is a completion of $\mc{A}$ and
Proposition~\ref{universal} shows that the natural morphism $\mc{U}(\mc{A}) \rightarrow G$ is surjective.

The completion $$\phi : \sqcup\mc{A} \rightarrow G$$ and, thus, the completion $$\hat\phi : \sqcup\mc{A}
\rightarrow \mc{U}(\mc{A})$$ is injective. Therefore the natural epimorphism $$\psi: \mc{U}(\mc{A}) \rightarrow
G$$ induces an isomorphism between the amalgam $\hat\phi(\sqcup\mc{A})$ inside $\mc{U}(\mc{A})$ and the amalgam
$\phi(\sqcup\mc{A})$ inside $G$. Hence the epimorphism $\psi : \mc{U}(\mc{A}) \rightarrow G$ induces a quotient map
between pregeometries $$\pi : \hat\G = ((\mc{U}(\mc{A})/G_{w_j^i}\times\{w_j^i\})_{i \in I, 1 \leq j \leq
t_i},*) \rightarrow ((G/G_{w_j^i}\times\{w_j^i\})_{i \in I, 1 \leq j \leq t_i},*).$$ The latter coset
pregeometry is isomorphic to $\G$ by the Reconstruction Theorem \ref{isomorphism theorem2}. Notice that
$\mc{U}(\mc{A})$ acts on $\G \cong ((G/G_{w_j^i}\times\{w_j^i\})_{i \in I, 1 \leq j \leq t_i},*)$ via
$$(gG_{w_j^i},w_j^i) \mapsto (\psi(u) gG_{w_j^i},w_j^i) \quad \mbox{ for } \quad u \in \mc{U}(\mc{A}).$$

We want to prove that this quotient map actually is a covering map. The pregeometry $\hat\G$ is connected by
Theorem \ref{char conn}, because $\mc{U}(\mc{A})$ is generated by $\hat{\phi}(\sqcup\mc{A})$. We establish the
isomorphism of the residues in a direct way. Fix an element $(G_{w_j^i},w_j^i)$ of $\G$, then, up to
isomorphism and with abuse of notation, $(G_{w_j^i},w_j^i)$ also denotes an arbitrary element of the inverse
image under $\pi$ of that element (we identify the elements $x$ and $\pi(x)$ when $x\in G_{w_{j^*}^{i^*}}$ for
some suitable $i^*$ and $j^*$).  In both $\G$ and $\hat\G$ the elements (not of type $i$) incident with
$(G_{w_j^i},w_j^i)$ can be written as $(uG_{w_{j'}^{i'}},w_{j'}^{i'})$, with $i'\neq i$, $1\leq j'\leq t_{i'}$
and $u\in G_{w_j^i}$. This already provides a natural bijection between the residues in $\G$ and $\hat\G$ of
$(G_{w_j^i},w_j^i)$. We now show that this bijection preserves incidence (in both directions). Clearly, if two
elements in $\hat\G$ are incident, then their projections in $\G$ are incident. So we are left to show that, if
$uG_{w_{j'}^{i'}}$ and $vG_{w_{j''}^{i''}}$ are non-disjoint in $G$, for $u,v\in G_{w^i_j}$, $i'\neq i\neq
i''$, $1\leq j'\leq t_{i'}$, $1\leq j''\leq t_{i''}$, then the ``same'' groups in $\mc{U}(\mc{A})$ are
non-disjoint (note that we identify $\widehat{\phi}(\sqcup\mc{A})$ with $\phi(\sqcup\mc{A})$ via the isomorphism induced by $\psi$ as above). By hypothesis (ii) and the fact that $w_j^i$, $w_{j'}^{i'}$ and $w_{j''}^{i''}$ are incident we
can assume that $u=v\in G_{w_j^i}$. But the definition of amalgams then implies that $uG_{w_{j'}^{i'}}$ and
$vG_{w_{j''}^{i''}}$ share the element $u=v$. Hence $\pi : \hat\G \rightarrow \G$ induces isomorphisms between
the residues of flags of rank one. So the map $\pi : \hat\G \rightarrow \G$ indeed is a covering of
pregeometries. Since $\G$ actually is a geometry the pregeometry $\hat\G$ is also a geometry.

Now we want to show that the covering $$\pi : \hat\G \rightarrow \G$$ induced by the canonical map
$\mc{U}(\mc{A}) \rightarrow G$ is universal. Denote the fundamental cover of $\G$ at some vertex $w^i_j$ of $W$
by $\G_0$ and let $$\phi : \G_0 \rightarrow \G$$ be the corresponding covering map. If $\hat w^i_j \in
\pi^{-1}(w^i_j)$, $\overline{w}^i_j \in \phi^{-1}(w^i_j)$, we will achieve the universality of $\pi$ by showing
that $\pi = \phi \circ \alpha$ for a unique isomorphism $$\alpha : \hat \G \rightarrow \G_0$$ with $\alpha(\hat
w^i_j) = \overline{w}^i_j$.
$$\xymatrix{
(\hat{\G}, \hat{w}_j^i) \ar[r]^\pi \ar[d]_\alpha & (\G, w_j^i) \\
(\G_0, \overline{w}_j^i) \ar[ur]_\phi & }$$ The simple connectivity of $\hat\G$ then is implied by the
universal property. For $g \in G_{w_j^i}$ define an automorphism $$\hat g^{(j,i)} : \G_0 \rightarrow \G_0 :
\phi^{-1}(\G,w^i_j) \ni \lbr \gamma \rbr \mapsto \lbr g(\gamma) \rbr.$$ The latter is also a homotopy class of
paths in $\G$ starting at $w^i_j$, because $g \in G_{w_j^i}$ stabilizes $w^i_j$. The fundamental cover $\G_0$
of $\G$ based at $w^i_j$ is isomorphic to the fundamental cover $\G_1$ of $\G$ based at some arbitrary
$w^{i'}_{j'} \in W$. Therefore we can define automorphisms on $\G_0$ using the automorphisms on $\G_1$ coming
from elements $g \in G_{w_{j'}^{i'}}$. To this end fix a maximal flag $V \subseteq W$ containing $w^i_j$. Let
$y \in V$ be incident to both $w^i_j$ and $w^{i'}_{j'}$ and for  $g \in G_{w_{j'}^{i'}}$ define an automorphism
$$\hat g^{(j',i')} : \G_0 \rightarrow \G_0 : (\lbr \gamma \rbr) \mapsto \lbr w^i_j, y, w^{i'}_{j'}, g(y),
g(\gamma) \rbr.$$ Since, for a different choice of $y' \in V$ incident to both $w^i_j$ and $w^{i'}_{j'}$, the
cycles $(y,y',w^i_j,y)$ and $(y,y',w^{i'}_{j'},y)$ are null-homotopic, the automorphism $\hat g^{(j',i')}$ does
not depend on the particular choice of $y \in V$. In particular, if $w^{i'}_{j'} \in V$, we can choose $y =
w^{i'}_{j'}$ or $y = w^i_j$.

Also, for incident $w^{i'}_{j'}$ and $w^{i''}_{j''}$, let $y$ be an element of $V$ incident to $w^i_j$,
$w^{i'}_{j'}$ and $w^{i''}_{j''}$. Since the cycles $(y,w^{i'}_{j'},w^{i''}_{j''},y)$ and
$(g(y),w^{i'}_{j'},w^{i''}_{j''},g(y))$ are null-homotopic, for $g \in G_{w_{j'}^{i'}} \cap G_{w_{j''}^{i''}}$
we have $$\lbr w^i_j, y, w^{i'}_{j'}, g(y), g(\gamma) \rbr = \lbr w^i_j, y, w^{i''}_{j''}, g(y), g(\gamma)
\rbr$$ and so $$\hat g^{(j',i')} = \hat g^{(j'',i'')}.$$

Hence  $$\hat{ } : \sqcup\mc{A} \rightarrow \hat G := \gen{\hat{\sqcup\mc{A}}} \leq \Aut\ \G_0$$ is a
completion map from $\mc{A}$ to $\hat G$. If $\hat g_1^{-1} \hat g_2$ acts trivially on $\G_0$, then
$g_1^{-1}g_2$ acts trivially on $\G$, thus $g_1=g_2$, as $G$ acts faithfully on $\G$. Therefore\ $\hat{ }$
embeds $\mc{A}$ in $\hat G$.

The geometry $\G_0$ together with the group $\hat G$ of automorphisms satisfies the hypothesis of the
Reconstruction Theorem \ref{isomorphism theorem2}, so  the geometry $\G_0$ is isomorphic to the coset
pregeometry $((\hat G/G_{w_{j}^{i}}\times\{w_j^i\})_{i \in I, 1 \leq j \leq t_i},*)$. The natural epimorphism
$\hat G \rightarrow G$ induces a covering map from $\G_0$ onto $\G$. Moreover, the natural epimorphism
$\mc{U}(\mc{A}) \rightarrow \hat G$ yields a quotient map $\hat\G \rightarrow \G_0$. Since $\G_0$ is universal
by Theorem \ref{1universal property} and therefore simply connected, this quotient map is a uniquely determined
isomorphism. Hence the covering $\pi : \hat\G \rightarrow \G$ is universal.

It remains to establish the structure of $\hat G \cong \mc{U}(\mc{A})$ to be of the form $\pi_1(\G).G$.
However, this is evident by Theorem \ref{1universal property}.
\end{proof}

We have now obtained a generalization of an important and widely used lemma by Tits.

\begin{corollary}[Tits' lemma] \label{tits} \label{Tits2} \label{tits2}
Let $\mc{G} = (X,*,\typ)$ be a geometry over $I$ of rank at least three and let $G$ be a group of automorphisms of $\G$. For each $i
\in I$ let $$w^i_1, \ldots, w^i_{t_i}$$ be $G$-orbit representatives of the elements of type $i$ of $\G$ such
that
\begin{enumerate}
\item $W := \bigcup_{i \in I} \lb w_1^i, \ldots, w_{t_i}^i \rb$ is a hall and, 
\item  if $V \subseteq W$ is a
chamber, the action of $G$ on the pregeometry over $\typ(V)$ consisting of all elements of the $G$-orbits $G.x$,
$x \in V$, is incidence-transitive and pennant-transitive.
\end{enumerate}
Let $\mc{A}(\G,G,W)$ be the amalgam of parabolics of $\G$ with respect to $G$ and $W$. The geometry $\G$ is
simply connected if and only if the canonical epimorphism $$\mc{U}(\mc{A}(\G,G,W)) \rightarrow G$$ is an
isomorphism. \pend
\end{corollary}

In practice it turns out that an amalgam may contain a lot of redundant information. Here is an explanation for that phenomenon.

\begin{theorem} \label{theorem444GHS}
Let $\mc{G} = (X,*,\typ)$ be a geometry over some finite set $I$ and let $G$ be a group of automorphisms of
$\G$. For each $i \in I$ let $$w^i_1, \ldots, w^i_{t_i}$$ be $G$-orbit representatives of the elements of type
$i$ of $\G$ such that
\begin{enumerate}
\item $W := \bigcup_{i \in I} \lb w_1^i, \ldots, w_{t_i}^i \rb$ is a hall and, \item  if $V \subseteq W$ is a
flag,
 the action of $G$ on the pregeometry over $\typ(V)$
consisting of all elements of the $G$-orbits $G.x$, $x \in V$, is
flag-transitive.
\end{enumerate}
Let $\mc{W} \subseteq 2^W$ be a shape, assume that for each flag $U \in 2^W \backslash \mc{W}$ the residue
$\G_U$ is simply connected, and let $\mc{A}(\G,G,W)$ and $\mc{A}_{\mc{W}}(\G,G,W)$ be the amalgam of maximal
parabolics respectively the amalgam of shape $\mc{W}$ of $\G$ with respect to $G$ and $W$. Then $$G =
\mc{U}(\mc{A}_{\mc{W}}(\G,G,W)).$$ In particular, if $\emptyset \not\in \mc{W}$, we have $$G =
\mc{U}(\mc{A}(\G,G,W)) = \mc{U}(\mc{A}_{\mc{W}}(\G,G,W)).$$
\end{theorem}

\begin{proof}
We will proceed by induction on the number of flags in the set $2^W \backslash \mc{W}$. If the set of flags
contained in $2^W \backslash \mc{W}$ is empty, then $\emptyset \in \mc{W}$, so the amalgam
$\mc{A}_{\mc{W}}(\G,G,W)$ contains the stabilizer in $G$ of the empty flag, i.e., $G$. Hence $G =
\mc{U}(\mc{A}_{\mc{W}}(\G,G,W))$. If there exists a flag in $2^W \backslash \mc{W}$, then the empty flag is
also contained in $2^W \backslash \mc{W}$, because by definition the shape $\mc{W}$ is closed under taking
superflags. Hence in that case $\G$ is simply connected and by Corollary \ref{Tits2} we have $G =
\mc{U}(\mc{A}(\G,G,W))$. We will now prove that $\mc{U}(\mc{A}(\G,G,W)) = \mc{U}(\mc{A}_{\mc{W}}(\G,G,W))$.

If the empty flag is the only flag contained in $2^W \backslash
\mc{W}$, then $\mc{A}(\G,G,W) = \mc{A}_{\mc{W}}(\G,G,W)$, so their
universal completions coincide. If there exists a nonempty flag in
$2^W \backslash \mc{W}$, then there also exists a (nonempty) flag
$U$ in $2^W \backslash \mc{W}$ such that $\mc{W}' := \lb U \rb
\cup \mc{W}$ is a shape. Then $\mc{A}_{\mc{W}'}(\G,G,W) =
\mc{A}_{\mc{W}}(\G,G,W) \cup G_U$. By connectivity of $\mc{G}_U$,
the group $G_U$ is a completion of the amalgam
$\mc{A}(\G_U,G_U,W_U)$, where $$W_U := W \cap \typ^{-1}(I
\backslash \typ(U)).$$ As $\G_U$ is simply connected, we even have
$$G_U = \mc{U}(\mc{A}(\G_U,G_U,W_U)).$$ Since $\mc{A}(\G_U,G_U,W_U)
\subseteq \mc{A}_{\mc{W}}(\G,G,W)$, we have
\begin{eqnarray*}
\mc{U}(\mc{A}_{\mc{W}}(\G,G,W)) & = & \mc{U}(\mc{A}_{\mc{W}}(\G,G,W) \cup \mc{U}(\mc{A}(\G_U,G_U,W_U))) \\ & =
&  \mc{U}(\mc{A}_{\mc{W}}(\G,G,W) \cup G_U) \\ & = & \mc{U}(\mc{A}_{\mc{W}'}(\G,G,W)).
\end{eqnarray*}
Hence, by induction, we have $\mc{U}(\mc{A}_{\mc{W}}(\G,G,W)) = \mc{U}(\mc{A}(\G,G,W))$, finishing the proof.
\end{proof}

The following corollary is a generalization of Theorem 8.2 of \cite{Gramlich/Hoffman/Nickel/Shpectorov:2005}.

\begin{corollary} \label{theorem44GHS}
Let $\mc{G} = (X,*,\typ)$ be a geometry over some finite set $I$, let $G$ be a group of automorphisms of $\G$,
for each $i \in I$ let $$w^i_1, \ldots, w^i_{t_i}$$ be $G$-orbit representatives of the elements of type $i$ of
$\G$ such that
\begin{enumerate}
\item $W := \bigcup_{i \in I} \lb w_1^i, \ldots, w_{t_i}^i \rb$ is a hall and,

\item if $V \subseteq W$ is a flag, the action of $G$ on the geometry $\typ(W)$ consisting of all elements of
the $G$-orbits $G.x$, $x \in V$, is flag-transitive.
\end{enumerate}
Let $k \leq |I|$, assume that all residues of rank greater or equal $k$ with respect to subsets of $W$ are
simply connected, and let $\mc{A}(\G,G,W)$ and $\mc{A}_{(k)}(\G,G,W)$ be the amalgam of maximal parabolics
respectively rank $k$ parabolics of $\G$ with respect to $G$ and $W$. Then $$G = \mc{U}(\mc{A}(\G,G,W)) =
\mc{U}(\mc{A}_{(k)}(\G,G,W)). $$ \pend
\end{corollary}

\section{Intransitive geometries: an example}

\subsection{Some standard techniques}
In this subsection, we collect some general results on simple
connectivity and null-homotopic cycles that have been established
in recent papers dealing with simple connectivity of
flag-transitive geometries.

\medskip
A \textbf{geometric cycle} in the geometry $\mc{G}$ is a cycle completely contained in the residue $\G_x$ of
some element $x$.

\begin{proposition}[Lemma 3.2 of \cite{Bennett/Shpectorov}] \label{geometric cycles}
Every geometric cycle is null-homotopic. \pend
\end{proposition}

\begin{corollary}[Lemma 3.3 of \cite{Bennett/Shpectorov}] \label{homotopy}
If two cycles based at the same element are obtained from one another by inserting or erasing a
geometric cycle then they are homotopic.\pend
\end{corollary}

\begin{definition}[Basic diagram] Let $\mc{G}$ be a geometry over
the set $I$. Let $i,j\in I$, then we define $i\sim j$ if there
exists a flag $f$ of cotype $\{i,j\}$ such that the residue of $f$
is a geometry containing two elements that are not incident. Then
the graph $(I,\sim)$ is called the \textbf{basic diagram of
$\mc{G}$}.

Let $\G$ be a geometry with basic diagram $$\node^{1}\stroke{}\node^{2}\quad\cdots,$$ i.e., the node $1$ has a
unique neighbor in the basic diagram of $\G$. Then the {\bf $1$-graph} (also called the {\bf collinearity
graph}) of $\G$ is the graph whose vertices are the elements of type $1$ (these elements will occasionally be
called \textbf{points}), where two such elements are adjacent if they are incident with a common element of
type $2$ (type 2 elements will sometimes be referred to as \textbf{lines}).
\end{definition}

\begin{definition}[Direct sum of pregeometries]
Let $\G = (X,*,\typ)$, $\G' = (X',*',\typ')$ be pregeometries over $I$ and $I'$. The {\bf direct sum} $$\G
\oplus \G'$$ is a pregeometry over $I \sqcup I'$
\begin{itemize}
\item whose element set is $X \sqcup X'$, \item whose type function is $\typ \cup \typ'$ and \item whose
incidence relation is the symmetric relation $*_\oplus$ with $*_\oplus|_{X \times X} = *$ and $*_\oplus|_{X'
\times X'} = *'$ and $*_\oplus|_{X \times X'} = X \times X'$, i.e., elements of $X$ are incident with elements
of $X'$.\end{itemize}
\end{definition}

For a geometry $\G$ over $I$, and for $J\subseteq I$, we denote by $_J\G$ the geometry over $J$, called {\bf truncation}, obtained from
$\G$ by deleting all elements of type not in $J$ (and retaining all other elements and incidences). For
instance, $_I\G=\G$ and $_\emptyset\G$ is the empty geometry.

\begin{lemma}[Lemma 5.1 of \cite{Gramlich/Hoffman/Shpectorov:2003}] \label{reduction22}
Let $\G$ be a geometry of rank $n\geq 3$ with basic diagram
$$\node^{1}\stroke{}\node^{2}\stroke{}\node_{}
\quad\cdots\quad\node_{}\stroke{}\node^{n}$$ and assume that for
each element $x$ of type $n$ the $1$-graph of $\G_x$ is connected.
Furthermore, suppose that if the residue $\G_x$ of some element
$x$ has a disconnected diagram falling into the two connected
components $\Delta_1$ and $\Delta_2$, then $\G_x$ is equal to the
direct sum $$_{\typ(\Delta_1)}\G_x \oplus {}_{\typ(\Delta_2)}\G_x.$$
Then every cycle of $\G$ based at some element of type $1$ or $2$
is homotopically equivalent to a cycle passing exclusively through
elements of type $1$ or $2$. \pend
\end{lemma}

\begin{lemma}[Lemma 7.2 of \cite{Gramlich/Hoffman/Shpectorov:2003}] \label{disconnected diagram}
Assume that $\G =\G_1\oplus \G_2$ with $\G_1$ connected of rank at
least two and $\G_2$ nonempty.  Then $\G$ is simply connected. \pend
\end{lemma}

Lemma~\ref{reduction22} allows to consider the \textbf{collinearity graph} (i.e., the graph with the points |
elements of type 1 | as the relation ``being incident with a common line'' | a line is an element of type 2 | as adjacency
relation) instead of the incidence graph when looking at cycles. Indeed, in order to prove simply connectivity
of certain geometries satisfying the conditions of the lemma, it is enough to deal with cycles consisting
merely of points and lines; such cycles can be considered as cycles in the collinearity graph (of half the
original length). Note that in the collinearity graph a triangle is not necessarily null-homotopic.

We will consider cycles in the collinearity graph from Lemma~\ref{orth3} below on.

\subsection{Generalities about orthogonal spaces}

Let $n \geq 1$ and let $V$ be an $(n+1)$-dimensional vector space
over some field $\mathbb{F}$ of characteristic distinct from $2$
endowed with some nondegenerate symmetric bilinear form
$f=(\cdot,\cdot)$. By $$\G_A^{\rm orth} = \G_A^{\rm
orth}(n,\mathbb{F},f)$$ we denote the pregeometry on the proper
subspaces of $V$ that are nondegenerate with respect to
$(\cdot,\cdot)$ with symmetrized containment as incidence and the
vector space dimension as the type.

\subsubsection*{Arbitrary fields of characteristic not two}

\begin{proposition}
The pregeometry $\G_A^{\rm orth}(n,\mathbb{F},f)$ is a geometry.
\end{proposition}

\begin{proof}
We have to prove that each flag can be embedded in a flag of cardinality $n$. To this end let $F = \left\{ x_1,
\ldots, x_t \right\}$ be a flag of $\G_A^{\rm orth}$. We can assume that the nondegenerate subspace $x_1$ of
$V$ has dimension one. Indeed, if it has not, then we can find a nondegenerate one-dimensional subspace $x_0$
of $x_1$ and study the flag $F' = F \cup \left\{ x_0 \right\}$ instead. Now observe that the residue of the
nondegenerate one-dimensional subspace $x_1$ is isomorphic to $\G_A^{\rm orth}(n-1,\mathbb{F},f')$ for some
induced form $f'$ via the map that sends an element $U$ of the residue of $x_1$ to $U \cap x_1^\perp$. Hence
induction applies.
\end{proof}

From now on, the notions of \emph{points} and \emph{lines} refer to points and lines, respectively, of the
geometry $\G_A^{\rm orth}$.

\begin{lemma}\label{pointline}
If $l$ is a line and $a$ is a point not on $l$, then there are at most two points of $\G_A^{\rm orth}$ on $l$
which are not collinear to $a$. In other words, if $\F$ is the field $\F_q$ of $q$ elements, there exist at
least $q-3$ points on $l$ collinear to $a$.
\end{lemma}

\begin{proof}
Let $U$ be the three-dimensional space $\<a,l\>$ and let $W=U\cap a^\perp$. The space $W$ has dimension two as $U$ has
dimension three. Hence there are at most two singular points on $W$ and, thus, there are at least $q-1$
nondegenerate lines in $U$ through $a$. The line $l$ has zero or two singular points, so at least $q-3$ of the
nondegenerate lines in $W$ through $a$ intersect $l$ is a nonsingular point.
\end{proof}

\begin{proposition}\label{diam}
Let $n \geq 3$ or $n=2$ and $| \F | \geq 5$. Then the collinearity graph of $\G_A^{\rm orth}(n,\F,f)$ has
diameter two.
\end{proposition}

\begin{proof}
If $n\ge 3$, then the dimension of the vector space $V$ is at least $4$. Fix two points
  $\la a\ra$ and $ \la b \ra$ which are not collinear, i.e.,
  the space $\la a,b\ra$ is singular with
  respect to $\form$. However $\la
  a,b \ra$ is a two-dimensional subspace of $ V$ which is not totally
  singular, because $(a,a)$ and $(b,b)$ are distinct from
  zero. Therefore the radical
  of $\la a,b \ra$ is a one-dimensional space. The dimension of $
  \la a, b \ra ^{\perp}$ is greater or equal to $2$.
  Consequently, the orthogonal complement of $\la a,b\ra$ contains a
  point, say
  $\la c \ra$. Consider the two two-dimensional subspaces $\la
  a,c\ra$ and $\la b,c\ra$. Since $\gen{a}$ and $\gen{b}$ are perpendicular to
  $\gen{c}$, both $\gen{a,c}$ and $\gen{b,c}$ are
  lines. The distance between $\gen{a}$ and $\gen{c}$
  is one and so is the distance between $\gen{c}$ and
  $\gen{b}$. Thus the distance between $\gen{a}$ and $\gen{b}$ is two.
Certainly $\G_A^{\rm orth}$ contains a pair of noncollinear points, so we have proved the claim for $n \geq 3$.

 If $n = 2$, let $\gen{a}$ and $\gen{b}$ be two arbitrary points in $V$. If
 the space $l=\gen{a,b}$ is a line then the distance between $\gen{a}$ and
 $\gen{b}$ is one. Otherwise pick a point $\gen{\tilde{a}}$ in
 $\gen{a}^{\perp}$. The space
 $\gen{a,\tilde{a}}$ is a line and the point $\gen{b}$ is not on
 $\gen{a,\tilde{a}}$.
 The point $\gen{b}$ is collinear with at least two points on
 $\gen{a,\tilde{a}}$ by Lemma \ref{pointline}. Pick one of these
 points, say the point
 $\gen{c}$. We have established that the
 distance between $\gen{a}$ and $\gen{b}$ is two.
\end{proof}

\begin{corollary}
Let $n \geq 2$ and $| \F | \geq 5$. Then $\G_A^{\rm orth}(n,\mathbb{F},f)$ is residually connected. \pend
\end{corollary}

It is shown in \cite{Altmann/Gramlich} that, if $n \geq 3$ and $\F$ not equal to $\F_3$ or $\F_5$, then the
geometry $\G_A^{\rm orth}(n,\mathbb{F},f)$ is simply connected. If the field $\F$ is quadratically closed, then
$\G_A^{\rm orth}(n,\mathbb{F},f)$ is flag-transitive and one can apply Corollary \ref{tits2} (Tits' lemma) to
obtain presentations of flag-transitive groups of automorphisms of that geometry, see \cite{Altmann/Gramlich}.
Also, in some cases like for real closed fields, it is possible to pass to suitable simply connected
flag-transitive parts of $\G_A^{\rm orth}(n,\mathbb{F},f)$ in order to obtain presentations of groups of
automorphisms.

\subsubsection*{Finite fields of characteristic not two}

The finite field case has been treated by Roberts in his PhD Thesis \cite{Rob}. He establishes the simple connectivity of certain suitably chosen flag-transitive parts of $\G_A^{\rm orth}(n,\mathbb{F},f)$ in case the rank of that geometry is at least four. Our approach allows us more freedom when choosing the subgeometry and enables us to deal with the case of rank three. 

\medskip We will be using standard terminology. In particular, each finite-dimensional vector space over some finite
field admits two isometry classes of nondegenerate quadratic forms, one called {\bf hyperbolic} (also {\bf
positive} or {\bf of plus type}), the other called {\bf elliptic} (also {\bf negative} or {\bf of minus type}).

Recall the following rules for determining the type of an orthogonal sum of nondegenerate orthogonal spaces
over a finite field $\mathbb{F}_q$. If $q \equiv 1\ \mod 4$ or if $q \equiv 3\ \mod 4$ and one of the involved
subspaces is even-dimensional, we have the following rule:
\begin{eqnarray*}
+ \oplus + & =  & +  , \\
+ \oplus - & = & -  , \\
- \oplus - & = & +  .
\end{eqnarray*}
If $q \equiv 3\ \mod 4$ and both of the involved subspaces are odd-dimensional, the following rule holds:
\begin{eqnarray*}
+ \oplus + & =  & -  , \\
+ \oplus - & = & +  , \\
- \oplus - & = & -  .
\end{eqnarray*}

The names ``hyperbolic'' and ``elliptic'' are a generalization of the classical usual incidence-theoretic
meaning: if a nondegenerate subspace of even dimension $2n\geq 2$ intersects the null-set of a quadratic form
in a quadric with Witt index $n$ or $n-1$, respectively, then the subspace is hyperbolic or elliptic,
respectively. We extend this as follows. If $q\equiv 1\ \mod 4$, and if a one-space takes square values or nonsquare values, respectively, with respect to the quadratic form, then this one-space is hyperbolic or elliptic,
respectively. Now these assignments of hyperbolic and elliptic, together with the above rules, determine the
plus/minus type of all nondegenerate subspaces (including the whole space and the zero space). Suppose now that
$q\equiv 3\ \mod 4$, and that $f$ is a non-degenerate quadratic form of a $d$-dimensional vector space over the
field of $q$ elements. Let $\Delta(f)$ be the discriminant of $f$ and put $e=\lfloor d/2\rfloor$ (hence $2e=d$
or $2e+1=d$). Then $f$ is positive if $\Delta(f)=(-1)^e\epsilon$, with $\epsilon$ a square if $d$ is even and
$\epsilon$ a non-square if $d$ is odd. Otherwise $f$ is negative.

\medskip
 The main tool for
our proof of simple connectivity is the following lemma. It is
clear that it would fail for transitive geometries as, roughly
speaking, one loses half the points when passing to a transitive
geometry.

\begin{lemma} \label{thelemmaorthfin}
Let $n \geq 2$, let $\F$ be a finite field of odd order $q$, let $p$ be a point of $\G_A^{\rm orth}(n,\F,f)$,
let $l$ be an elliptic line such that $\gen{p,l}$ is a nondegenerate plane, and let $m$ be a hyperbolic line
such that $\gen{p,m}$ is a nondegenerate plane. Then there exist at least $\frac{q-1}{2}$ elliptic lines
through $p$ intersecting $l$ in a point of $\G_A^{\rm orth}(n,\F,f)$ and at least $\frac{q-5}{2}$ hyperbolic
lines through $p$ intersecting $m$ in a point of $\G_A^{\rm orth}(n,\F,f)$.
\end{lemma}

\begin{proof}
Consider the two-dimensional nondegenerate space $p^\perp \cap
\gen{p,l}$. It contains $\frac{q+1}{2}$ or $\frac{q-1}{2}$ points
of positive type and $\frac{q+1}{2}$ or $\frac{q-1}{2}$ points of
$-$ type. Therefore, there exist at least $\frac{q-1}{2}$ elliptic
lines through $p$ intersecting $p^\perp \cap \gen{p,l}$ and, thus,
also $l$. The claim follows as all points on an elliptic line are
nondegenerate.

The number $\frac{q-5}{2} = \frac{q-1}{2} - 2$ of hyperbolic lines
through $p$ intersecting $m$ in a point of $\G_A^{\rm
orth}(n,\F,f)$ is obtained in exactly the same way plus the
observation that two of the hyperbolic lines through $p$ and $p^\perp
\cap \gen{p,m}$ could intersect $m$ in a singular point.
\end{proof}

\subsection{Positive form in dimension at least five}

Let $q$ be odd and let $V$ be a vector space over $\F_q$ of dimension $n+1$ at least five endowed with a
nondegenerate positive symmetric bilinear form $f$ and let $$\G_A^{\rm orth}(n,\F_q,f) = (X,*,\typ)$$ be the
pregeometry on all nondegenerate subspaces of $V$. Let $$W = \lb p, p', l, \pi, U, U_1, U_2, \ldots, U_t \rb$$
be a hall where $p$ is a positive point, $p'$ is a negative point, $l$ is a negative line, $\pi$ is a positive
or negative plane, $U$ is a positive four-dimensional subspace of $V$, and the $U_i$ are arbitrary
nondegenerate proper subspaces of $V$ of dimension at least three. Let $$(\G_A^{\rm orth})^W = (Y, *_{|Y \times
Y}, \typ_{|Y})$$ be a pregeometry with $$Y = \lb x \in X \mid \mbox{ there exists a $g \in \SO_{n+1}(\F_q,f)$
with $x \in g(W)$} \rb.$$

\begin{proposition} \label{everything holds}
The pregeometry $(\G_A^{\rm orth})^W$ is a geometry of rank
$|\typ(W)| \geq 3$ with linear diagram and a collinearity graph of
diameter two. Moreover, for each element $x$ of maximal type the
collinearity graph of the residue $(\G_A^{\rm orth})^W_x$ is
connected. Furthermore, if the residue $(\G_A^{\rm orth})^W_x$ of
some element $x$ has a disconnected diagram falling into the two
connected components $\Delta_1$ and $\Delta_2$, then $\G_x$ is
equal to the direct sum $$_{\typ(\Delta_1)}(\G_A^{\rm orth})^W_x
\oplus {}_{\typ(\Delta_2)}(\G_A^{\rm orth})^W_x.$$
\end{proposition}

\begin{proof}
In order to prove that $(\G_A^{\rm orth})^W$ is a geometry, it suffices to show that, for any two elements
$U,U'$ of type $i,i'$, respectively, with $i<i'$, there is an element $U^*$ of type $i+1$ incident with both
$U$ and $U'$. By taking the quotient projective space and corresponding quadratic form with respect to $U$ (or,
equivalently, by looking in $U^\perp$), we may assume that $U$ is the empty space. Hence $U^*$ is just some
positive or negative point in $U'$, which can always be found.

To prove the statement on the collinearity graph of $(\G_A^{\rm
orth})^W$ let $p$ and $p'$ be points of $(\G_A^{\rm orth})^W$.
Then there exists an elliptic line $l$ through $p'$ with
$\gen{p,l}$ nondegenerate. By Lemma \ref{thelemmaorthfin} there
exist $\frac{q-1}{2}$ elliptic lines through $p$ intersecting $l$
in a point of $(\G_A^{\rm orth})^W$. Since $q$ is odd, there
exists at least one, and the claim is proved. The same argument
implies that the collinearity graph of the residue of an element
$x$ of maximal type, which is at least four, is connected.
\end{proof}

The preceding proposition allows us to apply Lemma
\ref{reduction22}, so we can study the collinearity graph of
$(\G_A^{\rm orth})^W$ in order to establish the simple
connectivity of $(\G_A^{\rm orth})^W$.

In the following lemma we have to distinguish between $q \equiv 1\ \mod 4$ and $q \equiv 3\ \mod 4$. The latter
case is handled in parentheses.

\begin{lemma} \label{orth3}
Let $q > 7$. Then any triangle in the collinearity graph of $(\G_A^{\rm orth})^W$ is homotopically trivial.
\end{lemma}

\begin{proof}
Let $a$, $b$, $c$ denote the points of a triangle. If $\gen{a,b,c}$ is nondegenerate, then its polar
$\gen{a,b,c}^\perp$ contains a nondegenerate two-dimensional subspace of $V$ and, thus, points of positive type
and of negative type. Choosing a positive point $p$ of that line if $\gen{a,b,c}$ is positive [negative] and choosing a
negative point $p$ of that line if $\gen{a,b,c}$ is negative [positive], we obtain a positive space $\gen{a,b,c,p}$
containing the triangle $a$, $b$, $c$. Therefore that triangle is geometric, whence null-homotopic by
Proposition \ref{geometric cycles}.

Now suppose the triangle $a$, $b$, $c$ spans a degenerate space $\gen{a,b,c}$ with one-dimensional radical $x$.
Notice first that any line not passing through $x$ is elliptic. If $a$, $b$, $c$ are all of positive [negative] type
consider an arbitrary nondegenerate four-dimensional subspace of $V$ containing $\gen{a,b,c}$. That
four-dimensional space necessarily is of negative type, so its polar contains a negative point $p$. But
$\gen{a,p}$, $\gen{b,p}$, $\gen{c,p}$ then are elliptic lines and the three-dimensional spaces $\gen{a,b,p}$,
$\gen{b,c,p}$, $\gen{a,c,p}$ are nondegenerate, so the original triangle $a$, $b$, $c$ is null-homotopic. If
all of $a$, $b$, $c$ are negative [positive] points, then we can choose any positive [negative] point $p$ on the line $\gen{b,c}$
such that $\gen{a,p}$ does not contain $x$. Then $\gen{a,p}$ is an elliptic line and we have decomposed the
triangle $a$, $b$, $c$ into two triangles in which positive [negative] points occur. If $b$ and $c$ are of negative [positive] type
and $a$ is of positive [negative] type we can again choose any positive [negative] point $p$ on the line $\gen{b,c}$ such that
$\gen{a,p}$ does not contain $x$, decomposing the triangle $a$, $b$, $c$ into two triangles with one negative [positive]
point and two positive [negative] points.

We are left with the case of one negative [positive] point, say $a$, and two positive [negative] points, say $b$ and $c$. If neither
$b$ nor $c$ are perpendicular to $a$, we can choose the point $d$ on $\gen{b,c}$ perpendicular to $a$, which
is a positive point as it is perpendicular to the negative [positive] point on the nondegenerate, whence elliptic (i.e.\ negative) two-dimensional subspace $\gen{a,d}$. It remains to prove that the triangle $a$, $d$, $c$ is null-homotopic, the triangle $a$, $d$, $b$ being handled similarly. The space $\gen{a,d,c} = \gen{a,b,c}$ is contained in a four-dimensional nondegenerate
negative space which is in turn contained in a five-dimensional nondegenerate positive space $W$ (which may be
equal to $V$). The space $U := \gen{b,c}^\perp \cap W = \gen{d,c}^\perp \cap W$ is a three-dimensional negative space. As $d \perp a$
the space $\gen{a,U}$ equals $d^\perp \cap W$, which is a nondegenerate four-dimensional positive space.
Through $a$ there are $q+1$ tangent planes of $\gen{a,U}$. Moreover, in $U$ there are $q+1$ tangent lines. If
all tangent planes through $a$ would pass through a tangent line of $U$, then $U\subseteq a^\perp$, hence
$\gen{a,b,c}=U^\perp\cap W$. This would imply that $a$, $b$, $c$ are linearly dependent. So there exists a
nondegenerate plane of $\gen{a,U}$ through $a$ that intersects $U$ in a tangent line $L$ of $U$ with
corresponding tangent point $p_0$. Since $U$ is a negative space tangent lines of $U$ contain $q$ negative
points besides the radical. We have finished the proof, if we find a point $p$ among those $q$ points that spans an elliptic line
together with $a$ and nondegenerate three-dimensional spaces with $\gen{a,d}$ and $\gen{a,c}$, because then we can decompose the triangle $a$, $d$, $c$ into the three nondegenerate, whence null-homotopic triangles $a$, $d$, $p$ and $a$, $c$, $p$ and $c$, $d$, $p$. Since $d \perp
a$ and $d \perp p$, the space $\gen{a,d,p}$ is automatically nondegenerate if $\gen{a,p}$ is an elliptic line.
The space $\gen{a,c,p}$ has a Gram matrix of the form
$$\begin{pmatrix} * & * & \alpha \\ * & * & 0 \\ \alpha & 0 & \gamma \end{pmatrix}$$ with respect to the basis $a$,
$c$, $p$ for nonzero $\gamma$ and $\alpha$ both depending on $p$. If we choose a fixed point $p_1\neq p_0$ on
$L$ and write $p=\frac{1}{l}p_0+p_1$, for $l$ varying over the nonzero elements of the field, then we see that
$\gamma$ is constant and $\alpha$ runs through $q$ different elements of the field. Since the determinant of
the above matrix is quadratic in $\alpha$, we see that there are at most two choices of $p$ for which
$\gen{a,c,p}$ is degenerate. Hence there exist $k:=q - 2 - 2 - \frac{q-1}{2}$ points $p$ on a common elliptic
line with $a$. Indeed, there are $q$ negative points, two of which might give rise to a degenerate space
$\gen{a,c,p}$, two of which might give rise to a degenerate space $\gen{a,p}$ and $\frac{q-1}{2}$ of which
might span hyperbolic lines together with $a$. The number $k$ is positive since $q > 7$.
\end{proof}

\begin{lemma} \label{5424}
Let $q > 7$. Then any quadrangle of the collinearity graph of $(\G_A^{\rm orth})^W$ is homotopically trivial.
\end{lemma}

\begin{proof}
Let $a$, $b$, $c$, $d$ be a quadrangle and let $l := ab$ and $m := cd$. If $l$ and $m$ intersect in a point
$e$, then the quadrangle $a$, $b$, $c$, $d$ decomposes into two triangles $a$, $d$, $e$ and $b$, $c$, $e$.

Therefore we can assume $\gen{l,m}$ is four-dimensional. Our goal is to prove the claim that the point line
geometry consisting of the points of $l$ and $m$ and the elliptic lines in $\gen{l,m}$ intersecting $l$ and $m$
is connected. The fact that $a$, $b$, $c$, $d$ is null-homotopic then follows, as any path from $a$ to $b$ via
points on $l$ and $m$ and elliptic lines intersecting both $l$ and $m$ decomposes the quadrangle $a$, $b$, $c$,
$d$ into triangles. We have to consider the following five cases of possible structure for $\gen{l,m}$: (i)
two-dimensional radical, elliptic line as complement; (ii) two-dimensional radical, hyperbolic line as
complement; (iii) one-dimensional radical; (iv) nondegenerate negative space; (v) nondegenerate positive space.
In the first case any line not through the radical is elliptic and there is nothing to prove. The second case
cannot occur as the lines $l$ and $m$ are elliptic. In the third case let $x$ denote the radical of
$\gen{l,m}$. The planes $\gen{l,x}$ and $\gen{m,x}$ intersect in a line, $s$ say. Denote the intersection of
$l$ and $s$ by $y$ and the intersection of $m$ and $s$ by $z$. All lines in $\gen{l,x}$ through $z$ except $s$
are elliptic, whence $z$ is in the same connected component as any point on $l$ distinct from $y$. By symmetry,
$y$ is in the same connected component as any point on $m$ distinct from $z$. Now let $p$ be any point on $l$
distinct from $y$ and consider the plane $\gen{p,m}$. This plane is a complement in $\gen{l,m}$ of $x$, so it
is nondegenerate. By Lemma \ref{thelemmaorthfin} there exist $\frac{q-1}{2}$ elliptic lines through $p$ in
$\gen{p,m}$. This is at least two if $q$ is larger than three, so there exists an elliptic line through $p$
intersecting $m$ in a point distinct from $z$ and thus the claim follows. In case (iv) we can apply a similar
argument as above by using tangent planes of the elliptic quadric containing $l$ or $m$. In the fifth case the
space $\gen{l,m}$ is an object of the geometry $(\G_A^{\rm orth})^W$, so the quadrangle $a$, $b$, $c$, $d$ is
geometric and hence, by Lemma \ref{geometric cycles}, null-homotopic.
\end{proof}

\begin{lemma} \label{5425}
Suppose $q>7$. Any pentagon of the collinearity graph of $(\G_A^{\rm orth})^W$ is homotopically trivial.
\end{lemma}

\begin{proof}
Let $a$, $b$, $c$, $d$, $e$ be a pentagon and let $l := cd$. If $\gen{a,l}$ is nondegenerate, then there exist
$\frac{q-1}{2}$ elliptic lines through $a$ intersecting $l$, which is at least one, and if $\gen{a,l}$ is
degenerate, then there exist $q$ elliptic lines through $a$ intersecting $l$, as in $\gen{a,l}$ each complement
of the radical is an elliptic line. In both cases we have decomposed the pentagon $a$, $b$, $c$, $d$, $e$ into
two quadrangles.
\end{proof}

By Proposition \ref{everything holds}, the three lemmas we have proved yield the following theorem.

\begin{theorem} \label{finiteorth}
Let $q \geq 9$. Then the geometry $(\G_A^{\rm orth})^W$ is simply connected. \pend
\end{theorem}

\begin{theorem}
Let $q \geq 9$ be odd, let $n \geq 4$, let $V$ be an $(n+1)$-dimensional vector space over $\F_q$ endowed with
a nondegenerate positive symmetric bilinear form $f$. Let $\G = (\G_A^{\rm orth})^W$, let $G =
\SO_{n+1}(\F_q,f)$ and let $\mc{A} = \mc{A}(\G,G,W)$ be the amalgam of maximal parabolics of $(\G_A^{\rm
orth})^W$. Then $\mc{U}(\mc{A}) = \SO_{n+1}(\F_q,f)$.
\end{theorem}

\begin{proof}
This follows by Theorem \ref{finiteorth} and Corollary \ref{tits2}.
\end{proof}

\subsection{Negative form in dimension at least five}

Let $q$ be odd and let $V$ be a vector space over $\F_q$ of dimension $n+1$ at least five endowed with a
nondegenerate negative symmetric bilinear form $f$ and let $$\G_A^{\rm orth}(n,\F_q,f) = (X,*,\typ)$$ be the
pregeometry on all nondegenerate subspaces of $V$. Let $$W = \lb p, p', l, \pi, U, U_1, U_2, \ldots, U_t \rb$$
be a hall where $p$ is a positive point, $p'$ is a negative point, $l$ is a negative line, $\pi$ is a positive
or negative plane, $U$ is a positive four-dimensional subspace of $V$, and the $U_i$ are arbitrary
nondegenerate proper subspaces of $V$ of dimension at least three. Let $$(\G_A^{\rm orth})^W = (Y, *_{|Y \times
Y}, \typ_{|Y})$$ be a pregeometry with $$Y = \lb x \in X \mid \mbox{ there exists a $g \in \SO_{n+1}(\F_q,f)$
with $x \in g(W)$} \rb.$$

\begin{theorem} \label{finiteorth-}
Let $q \geq 9$. Then the geometry $(\G_A^{\rm orth})^W$ is simply connected.
\end{theorem}

\begin{proof}
The proof is almost the same as the proof of Theorem \ref{finiteorth}, i.e., it follows by versions of Lemmas
\ref{orth3}, \ref{5424} and \ref{5425}. The crucial step is finding a version of the proof of Lemma \ref{orth3}
that works. This, however, simply amounts to interchanging the words {\em positive} and {\em negative} in a
suitable way. The other two lemmas can be copied literally.
\end{proof}

\begin{theorem}
Let $q \geq 9$ be odd, let $n \geq 4$, let $V$ be an $(n+1)$-dimensional vector space over $\F_q$ endowed with
a nondegenerate negative symmetric bilinear form $f$. Let $\G = (\G_A^{\rm orth})^W$, let $G =
\SO_{n+1}(\F_q,f)$ and let $\mc{A} = \mc{A}(\G,G,W)$ be the amalgam of maximal parabolics of $(\G_A^{\rm
orth})^W$. Then $\mc{U}(\mc{A}) = \SO_{n+1}(\F_q,f)$. \pend
\end{theorem}

\subsection{Negative form in dimension four}

Let $q$ be odd and let $V$ be a vector space over $\F_q$ of dimension four endowed with a nondegenerate
negative symmetric bilinear form $f$ and let $$\G_A^{\rm orth}(3,\F_q,f) = (X,*,\typ)$$ be the pregeometry on
all nondegenerate subspaces of $V$. Let $$W = \lb p, p', l, \pi, \pi' \rb$$ be a hall where $p$ is a positive
point, $p'$ is a negative point, $l$ is a negative line, $\pi$ is a positive plane, and $\pi'$ is a negative
plane. Let $$(\G_A^{\rm orth})^W = (Y, *_{|Y \times Y}, \typ_{|Y})$$ be a pregeometry with $$Y = \lb x \in X
\mid \mbox{ there exists a $g \in \SO_{4}(\F_q,f)$ with $x \in g(W)$} \rb.$$

\begin{lemma}
Let $q \geq 9$. Then any triangle in the collinearity graph of $(\G_A^{\rm orth})^W$ is homotopically trivial.
\end{lemma}

\begin{proof}
Let $a,b,c$ be a triangle in a degenerate plane with one-dimensional radical $p$. There are two degenerate
planes through $bc$, namely $\gen{a,b,c}$ and some plane $\pi_{bc}$; likewise there are two degenerate planes
$\gen{a,b,c}$ and $\pi_{ac}$ through $ac$. The planes $\pi_{ac}$ and $\pi_{bc}$ meet in the line $L$ through
$c$. Since the tangent points of $\pi_{ac}$ and $\pi_{bc}$ are distinct, the line $L$ is an elliptic line (and
not a tangent line). There is a unique point $d_0\neq c$ on $L$ for which $\gen{a,b,d_0}$ is a degenerate
plane; there unique points $d_a$ and $d_b$ on $L$ for which $ad_a$ and $bd_b$ are degenerate. Since $q>7$, we
can pick a point $d\notin\{c,d_0,d_a,d_b\}$ on $L$. It follows that the plane $\gen{a,b,d}$ is nondegenerate.
Since there are precisely two degenerate planes through an elliptic line, there are two (not necessarily
distinct) points $c_a,c_b$ on the line $cp$ distinct from $c$ with the property that $\gen{ad,c_a}$ and
$\gen{b,d,c_b}$ are degenerate. Now, there are also $\frac{q-3}{2}$ points $c'$ on $cp$ distinct from $c$ (and
automatically distinct from $p$) such that $dc'$ is elliptic (noting that $\gen{c,d,p}$ cannot be degenerate as
$\gen{a,b,c}$ is the only degenerate plane containing $p$). As $\frac{q-3}{2}\geq 3$, we can choose $c'$
distinct from both $c_a$ and $c_b$. It is now clear that all triangles $\{a,b,d\}$, $\{a,c',d\}$ and
$\{b,c',d\}$ are contained in nondegenerate planes, and hence that $a$, $b$, $c'$ is null-homotopic. But
$\SO_{4}(\F_q,f)$ contains a group of order $q-1$ fixing $ab$ pointwise, fixing $p$ and acting transitively on
the points of $pc$ except for $p$ and the intersection $pc\cap ab$. So we conclude that also $a$, $b$, $c$ is
null-homotopic.
\end{proof}

\begin{theorem}
Let $q \geq 9$. Then $(\G_A^{\rm orth})^W$ is simply connected.
\end{theorem}

\begin{proof}
Case (iv) of Lemma \ref{5424} shows that any quadrangle of $(\G_A^{\rm orth})^W$ is null-homotopic and Lemma
\ref{5425} shows that any pentagon of $(\G_A^{\rm orth})^W$ is null-homotopic.
\end{proof}

\begin{theorem}
Let $q \geq 9$ be odd, let $V$ be a four-dimensional vector space over $\F_q$ endowed with a positive
nondegenerate form $f$. Let $\G = (\G_A^{\rm orth})^W$, let $G = \SO_{4}(\F_q,f)$ and let $\mc{A} =
\mc{A}(\G,G,W)$ be the amalgam of maximal parabolics of $(\G_A^{\rm orth})^W$. Then $\mc{U}(\mc{A}) =
\SO_{4}(\F_q,f)$. \pend
\end{theorem}

\subsection{Smaller amalgams}

\begin{theorem}
Let $q \geq 9$ be odd, let $n \geq 6$, let $V$ be an $(n+1)$-dimensional vector space over $\F_q$ endowed with
a nondegenerate positive symmetric bilinear form $f$. Assume that $W$ is a hall containing, besides a positive
and a negative point, a negative line, a positive or negative plane and a positive dimension four space, also a
positive and a negative hyperplane, a negative hyperline (which is a space of codimension two), a positive or
negative codimension three space and a positive codimension four space. Let $\G = (\G_A^{\rm orth})^W$, let $G
= \SO_{n+1}(\F_q,f)$ and let $\mc{A}_{n-2} = \mc{A}_{n-2}(\G,G,W)$ be the amalgam of rank $n-2$ parabolics of
$(\G_A^{\rm orth})^W$. Then
$$\mc{U}(\mc{A}_{n-2}) = \SO_{n+1}(\F_q,f).$$
\end{theorem}

\begin{proof}
In order to apply Corollary \ref{theorem44GHS}, we have to prove that the geometry itself and all
residues of flags of rank one are simply connected. Theorem \ref{finiteorth} implies that the geometry is simply connected. If the flag $x$ of rank one is neither a point nor a hyperplane of $(\G_A^{\rm
orth})^W$, then the simple connectivity of $(\G_A^{\rm orth})^W_x$ follows from
Lemma \ref{disconnected diagram}. If the flag $x$ is a positive hyperplane, then the simple connectivity of $(\G_A^{\rm orth})^W_x$ follows from Theorem \ref{finiteorth}. If the flag $x$ is a negative hyperplane, then the simple connectivity of $(\G_A^{\rm orth})^W_x$ follows from Theorem \ref{finiteorth-}.
In case $x$ being a point, we can dualize $(\G_A^{\rm orth})^W$ in order to reduce the situation to the case of $x$ being a hyperplane and, again, we can apply Theorem \ref{finiteorth}, resp.\ Theorem \ref{finiteorth-} to obtain simple connectivity.
\end{proof}

In principle, the theorem would also work for $n=5$, but then by assumption $W$ would have to contain a
negative line and a positive codimension four space, which would be a positive line. But this would contradict
the fact, that $W$ contains a positive and a negative point, because the connecting line between those two
points cannot be both positive and negative.

\begin{theorem}
Let $q \geq 9$ be odd, let $n \geq 4$, let $V$ be an $(n+1)$-dimensional vector space over $\F_q$ endowed with
a nondegenerate positive symmetric bilinear form $f$. Let
$$\G_A^{\rm orth}(n,\F_q,f) = (X,*,\typ)$$ be the pregeometry on
all nondegenerate subspaces of $V$. Let $$W = \lb p, p', l, \pi, U, U_1, U_2, \ldots, U_t \rb$$ be a hall where
$p$ is a positive point, $p'$ is a negative point, $l$ is a negative line, $\pi$ is a positive or negative
plane, $U$ is a positive four-dimensional subspace of $V$, and the $U_i$ are arbitrary nondegenerate proper
subspaces of $V$ of dimension at least three. Let $G = \SO_{n+1}(\F_q,f)$ and let $$(\G_A^{\rm orth})^W = (Y,
*_{|Y \times Y}, \typ_{|Y})$$ be a pregeometry  with $$Y = \lb x \in X \mid \mbox{ there exists a $g \in G$
with $x \in g(W)$} \rb.$$ Let $\mc{W} \subseteq 2^W$ be a shape containing $p$, $p'$, every flag of corank two,
and the flag consisting of all elements of type greater or equal four. Then $$G =
\mc{U}(\mc{A}_{\mc{W}}(\G,G,W)).$$
\end{theorem}

\begin{proof}
This follows from Corollary \ref{theorem44GHS}, Theorem \ref{finiteorth} and Lemmas \ref{tits2} and
\ref{disconnected diagram}.
\end{proof}

\section{Appendix: An intransitive geometry for $\mathsf{G_2}(3)$}

Here we present another application of our new theory. In \cite{Hoffman/Shpectorov} Hoffman and Shpectorov
study an amalgam of maximal subgroups of $\hat G = \Aut(\mathsf{G_2}(3))$ given by a certain choice of subgroups
\begin{eqnarray*}
\hat L & = &  2^3 \cdot \mathsf{L}_3(2) : 2 , \\
\hat N & = & 2^{1+4}.(S_3 \times S_3) , \\
M & = & \mathsf{G_2}(2) = \mathsf{U}_3(3):2
\end{eqnarray*}
which corresponds to an amalgam of subgroups of $G = G_2(3)$ given by
\begin{eqnarray*}
L & = &  \hat L \cap G = 2^3 \cdot \mathsf{L}_3(2) , \\
N & = & \hat N \cap G = 2^{1+4}.(3 \times 3).2 , \\
M & = & \mathsf{G_2}(2) = \mathsf{U}_3(3):2, \\
K & = & eMe^{-1} \quad \mbox{ for $e \in O_2(\hat L) \backslash O_2(L)$ }
\end{eqnarray*}
where $O_2(\hat L)$ denotes the largest normal subgroup of $\hat L$ that is a $2$-group.
The groups
\begin{eqnarray*}
\hat{G}_1 & = & \hat L, \\
\hat{G}_2 & = & \hat N, \\
\hat{G}_3 & = & M
\end{eqnarray*}
define a flag-transitive coset geometry $\G$ of rank three for $\hat G = \Aut(\mathsf{G_2}(3))$, which is simply
connected by \cite{Hoffman/Shpectorov}. The subgroup $G = \mathsf{G_2}(3)$ of $\hat G$ does not act flag-transitively on
$\G$. Nevertheless, the groups
\begin{eqnarray*}
G^{1,1} & = &  L,  \\
G^{1,2} & = & N, \\
G^{1,3} & = & M, \\
G^{2,3} & = & K
\end{eqnarray*}
define an intransitive coset geometry of rank three for $G = \mathsf{G_2}(3)$, which is isomorphic to $\G$ by
\cite{Hoffman/Shpectorov} and, hence, simply connected. Corollary \ref{Tits2} implies that $\hat G$ is the
universal completion of the amalgam given by $\hat L$, $\hat N$ and $M$ and their intersections as indicated in
Definition \ref{defrank} and that $G$ is the universal completion of the amalgam given by $L$, $N$, $M$ and $K$
and their intersections excluding $M \cap K$ as indicated in Definition \ref{defrank2}.

\noindent Authors' Addresses:

\medskip
\noindent \begin{tabular}{ll}Ralf Gramlich & Hendrik Van Maldeghem \\
TU Darmstadt & Ghent University\\
FB Mathematik / AG 5 & Pure Mathematics and Computer Algebra\\
Schlo\ss gartenstra\ss e 7 & Krijgslaan 281, S22\\
64289 Darmstadt & 9000 Gent\\
Germany & Belgium \\
{\tt gramlich@mathematik.tu-darmstadt.de}\phantom{lala}&{\tt hvm@cage.ugent.be}
\end{tabular}

\end{document}